\documentclass{article}
\usepackage{amsfonts,amsmath,amssymb,amsthm,txfonts}
\usepackage{extarrows}
\usepackage{graphicx}
\usepackage{geometry}
\usepackage{titlesec}

\allowdisplaybreaks

\begin{document}
\title{Heat kernel estimates on local and non-local Dirichlet spaces satisfying a weak chain condition}
\author{Guanhua Liu\footnote{The research is supported by SFB 1283 of the German Research Council.}}
\date{ }
\maketitle

\newtheorem{theorem}{Theorem}[section]
\newtheorem{proposition}[theorem]{Proposition}
\newtheorem{lemma}[theorem]{Lemma}
\newtheorem{remark}[theorem]{Remark}
\newtheorem{corollary}[theorem]{Corollary}
\numberwithin{equation}{section}

\begin{abstract}
In this paper, we focus on the heat kernel estimates for diffusions and jump processes on metric measure spaces satisfying a weak chain condition, where the length of a nearly shortest $\varepsilon$-chain between two points $x,y$ is comparable with a function of $d(x,y)$ and $\varepsilon$. For a diffusion, the best estimate is already given by Grigor'yan and Telcs, and we make it explicit in our particular case. For jump processes, especially those where the scale of the process is different with that of the jump kernel, we improve the results by Bae, Kang, Kim and Lee. Uniformity of the coefficients (or parameters) in the known estimates and metric transforms play the key role in our proof. We also show by examples how the weak chain condition is valid in practice.
\end{abstract}

\section{Introduction}

The heat kernel is an important tool in potential analysis on metric measure spaces, connecting Markovian processes and geometric analysis. Heat kernel estimates for diffusions provide much information on the metric measure structure and therefore constrain the range of regular jump processes. For example, on the Euclidean space $\mathbb{R}^d$, the Brownian motion is a diffusion corresponding to the classical Laplacian operator, which admits the following Gauss-Weierstrass function as its heat kernel:
$$p_t(x,y)=\frac{1}{(4\pi t)^{d/2}}\exp\left\{-\frac{|x-y|^2}{4t}\right\}.$$
Therefore, the behavior of this diffusion process is determined by $t$ and the ratio $|x-y|^2/t$. For this reason we call $t=r^2$ as its scale function. Correspondingly, a jump process with the jump kernel
$$J(x,y)=\frac{c}{|x-y|^{d+\beta}}$$
is regular if and only if $0<\beta<2$, and then its heat kernel satisfies the following stable-like estimate:
$$p_t(x,y)\asymp\frac{1}{t^{d/\beta}}\left(1+\frac{|x-y|^\beta}{t}\right)^{-(d+\beta)}\asymp\frac{1}{t^{d/\beta}}\wedge\frac{t}{|x-y|^{d+\beta}}.$$
Here $t=r^\beta$ is called the scale function of this jump process, and also of the jump kernel $J$. This way, the order of a regular jump process is lower than that of the diffusion. This can be explained by subordination theory.

Generally, on a metric measure space $(M,d,\mu)$, the scale function $\phi$ of a diffusion or a jump process is defined so that the heat kernel is, at least partly, determined by $t$ and $\phi(d(x,y))/t$ -- especially, for those places with $d(x,y)\le\phi^{-1}(t)$, which is called the near-diagonal region up to time $t$, the heat kernel $p_t(x,\cdot)$ is uniformly comparable with $\mu(B(x,\phi^{-1}(t)))^{-1}$. Generally $\phi$ may well not be a power function. However, under proper capacity conditions (for example, the Poincar\'e inequality and the cutoff Sobolev inequality), we can still obtain sub-Gaussian estimates (see \cite{ghl15}, or \cite{ab15} for the upper bound only)
\begin{equation}\label{sub-gauss}
p_t(x,y)\asymp\frac{1}{\mu(B(x,\phi^{-1}(t)))}\exp\left\{-\sup\limits_{\sigma>0}\left(c_1\frac{d(x,y)}{\sigma}-c_2\frac{t}{\phi(\sigma)}\right)\right\}
\end{equation}
and stable-like estimates (see \cite{ckw0}, or \cite{hl} for the upper bound only)
\begin{equation}\label{stb-lk}
p_t(x,y)\asymp\frac{1}{\mu(B(x,\phi^{-1}(t)))}\wedge\frac{t}{\mu(x,d(x,y))\phi(d(x,y))}.
\end{equation}

The parabolic Harnack inequality (for short, PHI) is slightly weaker than these estimates, actually equivalent to the upper estimate together with the near-diagonal lower estimate (see \cite{bbk} for diffusions and \cite{ckw} for jump processes). There are several ways to obtain the two-sided estimates (\ref{sub-gauss}) or (\ref{stb-lk}) from PHI, based on some extra conditions, for exmaple, the chain condition (for diffusions), and the assumption that the scale function $\phi$ of the process (usually consistent with the capacity conditions) is comparable with that of the jump kernel.

Some of these extra conditions are truly necessary for estimates (\ref{sub-gauss}) or (\ref{stb-lk}). The assumption on scale functions of jump processes is necessary by \cite{ckw0} clearly. As for a diffusion, \cite{mu} deduces a geometric condition from PHI, showing that (\ref{sub-gauss}) is impossible unless the chain condition holds.

However, there are many examples where such conditions fail to hold, but we still need to find proper heat kernel estimates. Significant progress has been made by \cite{gt} for diffusions without the chain condition, which we cite in Proposition \ref{p-phi-hk}; and by \cite{mi} for particular jump processes (to be exact, subordinated Brownian motions) on $\mathbb{R}^d$ with different scales for the jump kernel and the process. Generally, on a metric measure space equipped with a diffusion satisfying PHI with some scale function $\psi_c$, under the assumption
\begin{equation}\label{scl-int}
\int_0^1\frac{dt}{\psi_j\circ\psi_c^{-1}(t)}<\infty
\end{equation}
on the scale function $\psi_j$ of the jump kernel, we can establish the capacity conditions automatically with a new scale $\Phi$ (determined by $\psi_j$ and $\psi_c$) by subordination techniques. By \cite{lm}, (\ref{scl-int}) is an equivalent condition for a jump kernel with scale $\psi_j$ to define a regular jump process that satisfies PHI with some $\Phi$. \cite{bkkl} found two types of heat kernel estimates with both $\psi_j$ and $\Phi$, where the chain condition makes the different again: SHK (if the chain condition holds) and GHK (if not). The former is sharp, but the latter is not, though it is stronger than PHI.

This paper focuses on the case where $\Phi$ is not comparable with $\psi_j$ and the chain condition is not true. We establish a two-sided heat kernel estimate under a weaker condition on $\varepsilon$-chains (which we call the $\rho$-chain condition). It reduces to SHK if the (strong) chain condition holds, and is completely yielded by GHK and metric transforms. As a comparison, we also obtain a two-sided estimate on the diffusion under this $\rho$-chain condition, and in this way, the relation between these two estimates is obvious via subordination.

In our setting, $\infty$ stands for $+\infty$ only; for any $a\in\mathbb{R}$, $r>a$ and $r\ge a$ refer to $r\in(a,\infty)$ and $r\in[a,\infty)$ respectively (that is, $\infty$ is not allowed). Denote $\mathbb{R}_+:=[0,\infty)$, $a\wedge b:=\min\{a,b\}$ and $a\vee b:=\max\{a,b\}$ for any $a,b\in(-\infty,\infty]$. Let $\lceil a\rceil:=-[-a]=\min\{n\in\mathbb{Z}: n\ge a\}$ for any $a\in\mathbb{R}$. For an arbitrary non-negative function $f:A\to\mathbb{R}$, where $A$ is an arbitrary set, denote $\sup_{x\in\emptyset}f(x)=0$ and $\inf_{x\in\emptyset}f(x)=\infty$. We avoid ambiguity by not defining $\sup_{x\in\emptyset}f(x)$ for any $f$ that might take negative values. $C_0(M,d)$ is the collection of compactly supported functions on the metric space $(M,d)$. The value of a constant denoted by $C,C',c,c',c_1$ or $c_2$ may change from line to line, where generally $C,C'\ge 1$. Meanwhile, if a constant is denoted by $C_i,C'_i,c_i$ (except $c_1$ and $c_2$) or $c'_i$, its value is fixed once it appears. Whatever class $\mathfrak{K}$ is in, and whatever the subscript $i$ is, a constant is denoted as $C(\mathfrak{K}),C'(\mathfrak{K})$ or $c(\mathfrak{K})$ to indicate that it depends on $\mathfrak{K}$, but the dependence changes at every appearance; while the dependence of a constant denoted as $C_i(\mathfrak{K}),C'_i(\mathfrak{K})$ or $c_i(\mathfrak{K})$ on $\mathfrak{K}$ is fixed once it appears. Given constants $C,c_1,c_2>0$, a formula $h\asymp F(c_\bullet)$ refers to $C^{-1}F(c_1)\le h\le CF(c_2)$; while without showing $C$, $a\asymp b$ means $C^{-1}b\le a\le Cb$ with some constant $C$ that only (implicitly) appears in this line.

Let $(M,d,\mu)$ be a metric measure space, where $(M,d)$ is separable, locally compact and complete (so that every metric ball is precompact) and $\mu$ is a Radon measure with full support. We say the condition $\mathrm{RVD}$ holds, if there exist $C=C_{RVD}(M,d,\mu)$ and $\alpha_1=\alpha_1(M,d,\mu)$ such that for all $x\in M$ and $0<r\le R<\overline{R}:=diam(M,d)$,
$$\frac{V(x,R)}{V(x,r)}\ge C^{-1}\left(\frac{R}{r}\right)^{\alpha_1};$$
similarly, we say the condition $\mathrm{VD}$ holds, if there exist $C=C_{VD}(M,d,\mu)$ and $\alpha_2=\alpha_2(M,d,\mu)$ such that for all $x\in M$ and $0<r<R<\infty$,
$$\frac{V(x,R)}{V(x,r)}\le C\left(\frac{R}{r}\right)^{\alpha_2}.$$

Given two arbitrary points $x,y\in M$, for any $\varepsilon>0$, a sequence of points $\xi_0=x,\xi_1,\cdots,\xi_N,\xi_{N+1}=y$ is called an $\varepsilon$-chain from $x$ to $y$, if $\xi_i\in M$ and $d(\xi_i,\xi_{i+1})<\varepsilon$ for all $0\le i\le N$. Let
$$d_\varepsilon(x,y)=\inf\left\{\sum\limits_{i=0}^Nd(\xi_i,\xi_{i+1}):\{\xi_i\}_{i=0}^{N+1}\mbox{ is an }\varepsilon\mbox{-chain from }x\mbox{ to }y\right\}.$$
In particular, $d_\varepsilon(x,y)=\infty$ if there does not exist any $\varepsilon$-chain from $x$ to $y$. We say that $(M,d)$ satisfies the chain condition, if there exists a constant $C>0$ such that for all $x,y\in M$ and all $\varepsilon>0$, $d_\varepsilon(x,y)\le Cd(x,y)$.

$\phi:\mathbb{R}_+\to\mathbb{R}_+$ is called a scale function, if it is continuous, monotonically increasing, mapping 0 to 0, $\infty$ to $\infty$, and there exist constants $C(\phi)\ge 1$, $\beta_1(\phi),\beta_2(\phi)>0$ such that for any $0<s\le t<\infty$,
\begin{equation}\label{scale}
C(\phi)^{-1}\left(\frac{t}{s}\right)^{\beta_1(\phi)}\le\frac{\phi(t)}{\phi(s)}\le C(\phi)\left(\frac{t}{s}\right)^{\beta_2(\phi)}.
\end{equation}
Given $R\in(0,\infty]$, we take $\beta_1^{(R)}(\phi)\ge\beta_1(\phi)$ and $C(\phi,R)\ge C(\phi)$ satisfying
\begin{equation}\label{scale'}
\frac{\phi(t)}{\phi(s)}\ge C(\phi,R)^{-1}\left(\frac{t}{s}\right)^{\beta_1^{(R)}(\phi)},\quad\mbox{for all }0<s\le t<R.
\end{equation}
in a way such that $C(\phi,R)\le C'(\phi)C(\phi)$ for all $R\in(0,\infty]$ with some $C'(\phi)>0$, $\beta_1^{(R)}\le\beta_1^{(r)}$ if $R\ge r>0$, and $\beta_1^{(\infty)}(\phi)=\beta_1(\phi)$. By \cite[Lemma A.2]{bkkl}, if (\ref{scale'}) holds, then for all $0<s<t<\phi(R)$,
\begin{equation}\label{scale''}
\frac{\phi^{-1}(t)}{\phi^{-1}(s)}\le\left(C(\phi)\frac{t}{s}\right)^{1/\beta_1^{(R)}}.
\end{equation}
Throughout this paper, if a scale function is denoted by $\rho$, we assume that $\gamma_i:=\beta_i(\rho)\ge 1$. We say the $\rho$-chain condition holds, if there exists a constant $C=C(d,\rho)>0$ such that for all $x,y\in M$ and all $0<\varepsilon<d(x,y)$,
\begin{equation}\label{rho-chain}
C^{-1}\frac{\rho(d(x,y))}{\rho(\varepsilon)}\le\frac{d_\varepsilon(x,y)}{\varepsilon}\le C\frac{\rho(d(x,y))}{\rho(\varepsilon)}.
\end{equation}
By taking $\rho_1(r)=r$, clearly the chain condition is equivalent to the $\rho_1$-chain condition. In addition, if $(M,d)$ is connected, then the $\rho$-chain condition can be expressed equivalently as follows: there exists a constant $C>0$ such that for all $x,y\in M$ and all $0<\varepsilon<d(x,y)$, there exists an $\varepsilon$-chain $\{\xi_i\}_{i=0}^{N+1}$, with $N\le C\rho(d(x,y))/\rho(\varepsilon)$, from $x$ to $y$, but there does not exist such chains with $N\le C^{-1}\rho(d(x,y))/\rho(\varepsilon)$.

A function $J:M\times M\setminus\mathrm{diag}\to\mathbb{R}_+$ is called a jump kernel on $M$ if it is symmetric and measurable with respect to $\mu\otimes\mu$. As usual, we set $J(x,x)=0$ for all $x\in M$ so that a jump kernel is defined on $M\times M$. We say the condition $\mathrm{J}(\psi_j)$ holds, if $\psi_j$ is a scale function, and there exists a constant $C=C_J(J,\psi_j)$ such that for every two different points $x,y\in M$,
$$\frac{C^{-1}}{V(x,d(x,y))\psi_j(d(x,y))}\le J(x,y)\le \frac{C}{V(x,d(x,y))\psi_j(d(x,y))}.$$
We say $\mathrm{J}_\le(\psi_j)$ holds if the second inequality holds.

Let $(\mathcal{E},\mathcal{F})$ be a (symmetric) Dirichlet form on $(M,d,\mu)$, that is, $\mathcal{F}$ is a dense subspace of $L^2(M,\mu)$, $\mathcal{E}$ is a non-negative definite quadratic form on $\mathcal{F}$, $\mathcal{F}$ is a Hilbert space equipped with the inner product $\mathcal{E}_1(\cdot):=\mathcal{E}(\cdot,\cdot)+\|\cdot\|^2_{L^2(\mu)}$, and the Markovian property $\mathcal{E}(u_+\wedge 1,u_+\wedge 1)\le\mathcal{E}(u,u)$ holds. We only consider Dirichlet forms that are regular, i.e., $\mathcal{F}\cap C_0(M,d)$ is dense both in $\mathcal{F}$ (under $\mathcal{E}_1$) and in $C_0(M,d)$ (under the uniform norm). In this case, we also say that $(M,d,\mu,\mathcal{E},\mathcal{F})$ is a regular Dirichlet space. Let $\{P_t\}_{t>0}$ be its heat semigroup and $\Gamma$ be its energy form, that is, for all $u,v\in\mathcal{F}$, $\varphi\in\mathcal{F}\cap C_0(M,d)$, $f\in L^2(M,\mu)$ and $t>0$, we have $P_tf\in\mathcal{F}$ and $\Gamma(v)$ is a positive Radon measure on $M$ such that
$$\frac{d}{dt}\int_M\varphi P_tfd\mu+\mathcal{E}(\varphi,P_tf)=0,\quad\int_Mud\Gamma(v)=\mathcal{E}(uv,v)-\frac{1}{2}\mathcal{E}(u,v^2).$$
Let the Beurling-Deny decomposition of $(\mathcal{E},\mathcal{F})$ (cf. \cite[Theorem 3.2.1]{fot}) be
\begin{equation}\label{bd-fml}
\mathcal{E}(u,v)=\mathcal{E}^{(L)}(u,v)+\int_{M\times M}(u(x)-u(y))(v(x)-v(y))dj(x,y)+\int_Mu(x)v(x)dk(x),
\end{equation}
where $\mathcal{E}^{(L)}$ is strongly local (i.e., $\mathcal{E}^{(L)}(u,v)=0$ whenever $u,v\in\mathcal{F}$ and $v$ is a constant on an open neighborhood of supp$(u)$), and $j,k$ are positive Radon measures respectively on $M\times M\setminus\mathrm{diag}$ and on $M$. In this way, $\Gamma$ is correspondingly decomposed into $\Gamma^{(L)}$, $\Gamma^{(J)}$ and $\Gamma^{(K)}$.

If $j,k=0$, this Dirichlet form is strongly local, and the corresponding Hunt process $X_t$ (such that, starting from $x$, the transition probability at time $t$ onto $\Omega$ is $P_t1_\Omega(x)$) is then called a diffusion; while if $\mathcal{E}^{(L)},k=0$, this Dirichlet form is purely jump-type, and the corresponding Hunt process $X_t$ is called a jump process. In either case, let $\mathcal{F}':=\{u=u_0+a:\ u\in\mathcal{F}\mbox{ and }a\mbox{ is a constant}\}.$ The integration $(u,v)\mapsto\int_Muvdk$ is called the killing part of $(\mathcal{E},\mathcal{F})$, which corresponds to a killed process.

For any open subset $\Omega\subset M$, let $P_t^\Omega$ be the heat semigroup restricted on $\Omega$ (that is, a linear operator on $L^2(\Omega,\mu)$ such that the heat equation above holds for all $\varphi\in\mathcal{F}\cap C_0(\Omega,d)$ and $f\in L^2(\Omega,\mu)$, with $P_t^\Omega$ replacing $P_t$), and $G^\Omega f:=\int_0^\infty P_t^\Omega fdt$ be the Green operator on $\Omega$. If $dj(x,y)=J(x,dy)d\mu(x)$, then let
$$d\Gamma^\Omega(u)(x):=d\Gamma^{(L)}(u)(x)+\left\{\int_M(u(x)-u(y))^2J(x,dy)\right\}d\mu(x)+\frac{1}{2}u^2(x)dk(x),\mbox{ for all }x\in\Omega,$$
so that $d\Gamma^\Omega(u)\uparrow d\Gamma(u)$ as $\Omega$ exhausts $M$.

Based on these notions, we recall some standard conditions describing a Dirichlet space $(M,d,\mu,\mathcal{E},\mathcal{F})$ without killing part (where $\phi$ is a scale function):

\begin{itemize}
\item $\mathrm{PHI}(\phi)$ (the parabolic Harnack inequality): there exist constants $0<C_1<C_2<C_3<C_4<\infty$, $C=C_{PHI}(\mathcal{E},\mathcal{F},\phi)\ge 1$ and $c=c_{PHI}(\mathcal{E},\mathcal{F},\phi)\in(0,1)$ such that for all $0<t<(C_4^{-1}\wedge 1)\phi(\overline{R})$, $x\in M$ and all $u:\mathbb{R}_+\times M\to\mathbb{R}_+$ that is caloric on $(0,t]\times B(x,\phi^{-1}(t))$ (i.e. $\frac{d}{ds}(u(s,\cdot),\varphi)=-\mathcal{E}(u(s,\cdot),\varphi)$ for all $0<s\le t$ and $\varphi\in\mathcal{F}\cap C_0(B(x,\phi^{-1}(t)),d)$),
$$\sup\limits_{(C_1t,C_2t]\times B(x,c\phi^{-1}(t))}u\le C\inf\limits_{(C_3t,C_4t]\times B(x,c\phi^{-1}(t))}u;$$

\item $\mathrm{E}(\phi)$ (the exit time estimate): there exist constants $C=C_E(\mathcal{E},\mathcal{F},\phi)\ge 1$ and $\kappa=\kappa_E(\mathcal{E},\mathcal{F},\phi)\in(0,1]$ such that for all $x\in M$ and $0<r<\kappa\overline{R}$,
$$C^{-1}\phi(r)1_{B(x,r/4)}\le G^{B(x,r)}1_{B(x,r)}\le C\phi(r);$$

\item $\mathrm{DUE}(\phi)$ (the on-diagonal upper estimate): $P_t$ admits a heat kernel $p_t:M\times M\to\mathbb{R}_+$, and there exists a constant $C=C_{DUE}(\mathcal{E},\mathcal{F},\phi)\ge 1$ such that for all $x,y\in M$ and $0<t<\phi(\overline{R})$,
$$p_t(x,y)\le\frac{C}{V(x,\phi^{-1}(t))};$$

\item $\mathrm{NLE}(\phi)$ (the near-diagonal lower estimate): the heat kernel exists, and there exist constants $\eta\in(0,1)$ and $C=C_{NLE}(\mathcal{E},\mathcal{F},\phi)\ge 1$ such that for all $0<t<\phi(\overline{R})$ and $x,y\in M$ such that $d(x,y)\le\eta\phi^{-1}(t)$,
$$p_t(x,y)\ge\frac{C^{-1}}{V(x,\phi^{-1}(t))};$$

\item $\mathrm{LLE}(\phi)$ (the local lower estimate, also denoted as $\mathrm{NDL}(\phi)$ in some literature): the heat kernel exists everywhere, and there exist constants $\eta\in(0,1)$ and $C=C_{LLE}(\mathcal{E},\mathcal{F},\phi)\ge 1$ such that for all $0<t<\phi(\overline{R})$, every ball $B=B(x_0,R)$ with $R\ge\eta^{-1}\phi^{-1}(t)$, and all $x,y\in B(x_0,\eta\phi^{-1}(t))$,
$$p_t^B(x,y)\ge\frac{C^{-1}}{V(x_0,\phi^{-1}(t))};$$

\item $\mathrm{UJS}$ (the upper jumping smoothness): there exists $C=C_{UJS}(J)\ge 1$ such that for all $x,y\in M$ and all $0<r\le d(x,y)/2$,
$$J(x,y)\le C\fint_{B(x,r)}J(z,y)d\mu(z);$$

\item $\mathrm{ABB}(\phi)$ (the Andres-Barlow-Bass condition, also denoted as $\mathrm{CSJ}(\phi)$ in some literature): there exist $\zeta=\zeta_{ABB}(\mathcal{E},\mathcal{F},\phi)\ge 0$ and $C=C_{ABB}(\mathcal{E},\mathcal{F},\phi)>0$ such that for all concentric balls $B_i=B(x,r_i)$ ($i=1,2,3$, where $r_1<r_2<r_3$), there exists a cutoff function $\varphi\in C_0(B_2,d)$ for $B_1$ (that is, $\varphi|_{B_1}=1$ and $0\le\varphi\le 1$ on $B_2$) such that for all $u\in\mathcal{F}'\cap L^\infty(M,\mu)$,
$$\int_{B_3}u^2d\Gamma^{B_3}(\varphi)\le\zeta\int_{B_3}d\Gamma^{B_3}(u)+\frac{C}{\phi(r_2-r_1)}\int_{B_3}u^2d\mu;$$

\item $\mathrm{PI}(\phi)$ (the Poincar\'e inequality): there exist constants $C=C_{PI}(\mathcal{E},\mathcal{F},\phi),K=K_{PI}(\mathcal{E},\mathcal{F},\phi)\ge 1$ such that for all $x\in M$, $0<r<K^{-1}\overline{R}$ and $f\in\mathcal{F}$,
$$\int_{B(x,r)}\left(f(\xi)-\fint_{B(x,r)}f(y)d\mu(y)\right)^2d\mu(\xi)\le C\phi(r)\int_{B(x,Kr)}d\Gamma^{B(x,Kr)}(f);$$

\item $\mathrm{C}$ (the conservativeness, or the stochastic completeness): $P_t1=1$ for all $t>0$.
\end{itemize}

Here are some remarks on these conditions:

\begin{remark}\label{release-T}
For any $K\ge 1$, $\mathrm{PHI}(\phi)$ implies $\mathrm{PHI}(K\phi)$ by iteration, while $\mathrm{DUE}(\phi)$ implies $\mathrm{DUE}(K\phi)$ by the contractiveness of $P_t$ and the semigroup property under condition $\mathrm{VD}$, with $C_1,C_2,C_3,C_4$ invariant and
\begin{align*}
C_{PHI}(\mathcal{E},\mathcal{F},K\phi)=&C_{PHI}(\mathcal{E},\mathcal{F},\phi)^{\frac{C_4-C_1}{C_3-C_1}K+1},\quad c_{PHI}(\mathcal{E},\mathcal{F},K\phi)=K^{-1}c_{PHI}(\mathcal{E},\mathcal{F},\phi),\\
C_{DUE}(\mathcal{E},\mathcal{F},K\phi)=&C_{DUE}(\mathcal{E},\mathcal{F},\phi)C_{VD}C(\phi)^{1/\beta_1}K^{\alpha_2/\beta_1}.
\end{align*}
\end{remark}

The condition $\mathrm{E}(\phi)$ is called the exit time estimate because $G^B1_B(x)=\mathbb{E}^x[\tau_B]$, where $\tau_B:=\inf\{t>0: X_t\notin B\}$ is the stopping time when the Hunt process $X_t$ exits $B$ for the first time.

If $J$ is a jump kernel on $M$, and
$$\mathcal{E}(u,v):=\int_M\int_M(u(x)-u(y))^2J(x,y)d\mu(y)d\mu(x),\quad\mathcal{F}:=\{u\in L^2(M,\mu):\mathcal{E}(u,u)<\infty\}$$
defines a regular Dirichlet form, then the coefficients (or parameters) above are also expressed as $C_E(J,\phi)$, $\kappa_E(J,\phi)$, etc. Note that if two jump kernels $J_1\asymp J_2$, then $\mathrm{PI}(\phi)$ (or $\mathrm{ABB}(\phi)$) holds for $J_1$ if and only if it holds for $J_2$. This can be easily seen from the formula
$$d\Gamma^\Omega(u)(x)=1_\Omega(x)\left(\int_\Omega(u(x)-u(y))^2J(x,y)d\mu(y)\right)d\mu(x).$$

Additionally, denoted the upper estimates of type (\ref{sub-gauss}) and (\ref{stb-lk}) respectively as $\mathrm{UHK}_c(\phi)$ and $\mathrm{UHK}_j(\phi)$. Clearly $\mathrm{DUE}(\phi)$ is contained in either upper estimate.

We never use the definition of $\mathrm{PHI}(\phi)$ throughout this paper, but instead, refer to its equivalent conditions:

\begin{proposition}\label{PHI}
Under conditions $\mathrm{VD}$ and $\mathrm{RVD}$, if a regular Dirichlet form $(\mathcal{E},\mathcal{F})$ is strongly local, then
$$\mathrm{PHI}(\phi)\Leftrightarrow\mathrm{UHK}_c(\phi)+\mathrm{NLE}(\phi)\Leftrightarrow\mathrm{EHI}+\mathrm{E}(\phi)\Rightarrow\mathrm{PI}(\phi)+\mathrm{ABB}(\phi)\Rightarrow\mathrm{cap}(\phi)+\mathrm{C};$$
while if it is purely jump-type, then
$$\mathrm{PHI}(\phi)\Leftrightarrow\mathrm{UHK}_j(\phi)+\mathrm{NLE}(\phi)+\mathrm{UJS}\Leftrightarrow\mathrm{LLE}(\phi)+\mathrm{UJS}\Leftrightarrow\mathrm{J}_\le(\phi)+\mathrm{UJS}+\mathrm{PI}(\phi)+\mathrm{ABB}(\phi)\Rightarrow\mathrm{E}(\phi)+\mathrm{C}.$$
\end{proposition}

\begin{proof}
If $(M,d)$ is unbounded, the equivalences for the strongly local case come from \cite[Theorem 2.15]{bbk} and \cite[Theorem 1.2]{ghl15}; $\mathrm{cap}(\phi)$ comes from \cite[Theorem 1.1]{ghl15}. The equivalences and $\mathrm{E}(\phi)$ come from \cite[Theorem 1.18]{ckw}. See \cite[Remark 1.19]{ckw} how these implications fit the bounded spaces. Each case $\mathrm{C}$ comes from \cite[Theorem 1.3]{hl} (where $\mathrm{FK}(\phi)$ comes from $\mathrm{DUE}(\phi)$ by \cite[Lemma 6.8]{hl}).
\end{proof}

Here the conditions $\mathrm{cap}$, $\mathrm{FK}(\phi)$ and $\mathrm{EHI}$ appear only literally and keep consistent with the cited context. Please check the citations for their definitions.

The paper aims to find the best heat kernel estimates for local and non-local (especially the latter) regular Dirichlet forms on metric measure spaces satisfying some $\rho$-chain condition. The following notions on heat kernel estimates come from \cite{bkkl}. Given scale functions $\psi_c$ and $\psi_j$ satisfying (\ref{scl-int}), define
\begin{equation}\label{phi}
\Phi(r):=\frac{\psi_c(r)}{\int_0^r\frac{\psi_c(s)ds}{s\psi_j(s)}}.
\end{equation}
We say the condition $\mathrm{SHK}(\Phi,\psi_j;T)$ holds, if there exist $C,c_1,c_2>0$ such that for all $x,y\in M$ and $0<t<T\wedge\Phi(\overline{R})$,
\begin{equation}\label{shk}
p_t(x,y)\asymp\frac{1}{V(x,\Phi^{-1}(t))}\wedge\left\{\frac{t}{V(x,d(x,y))\psi_j(d(x,y))}+\frac{1}{V(x,\Phi^{-1}(t))}\exp\left(-c_\bullet\sup\limits_{\sigma>0}\left[\frac{d(x,y)}{\sigma}-\frac{t}{\Phi(\sigma)}\right]\right)\right\}.
\end{equation}
We say $\mathrm{UHK}(\Phi,\psi_j;T)$ holds, if the upper bound in (\ref{shk}) holds for all $x,y\in M$ and $0<t<T\wedge\Phi(\overline{R})$. We say $\mathrm{G^+HK}(\Phi,\psi_j;T)$ holds, if $\mathrm{UHK}(\Phi,\psi_j;T)$ holds and there exist $C,\eta>0$ such that for all $x,y\in M$ and $0<t<T\wedge\Phi(\overline{R})$,
$$p_t(x,y)\ge\begin{cases}\frac{C^{-1}}{V(x,\Phi^{-1}(t))},&\mbox{if }d(x,y)\le\eta\Phi^{-1}(t);\\
\frac{C^{-1}t}{V(x,d(x,y))\psi_j(d(x,y))},&\mbox{if }d(x,y)>\eta\Phi^{-1}(t).
\end{cases}$$
When $T\in[\Phi(\overline{R}),\infty]$, we omit this parameter and write $\mathrm{SHK}(\Phi,\psi_j)$, $\mathrm{UHK}(\Phi,\psi_j)$ and $\mathrm{G^+HK}(\Phi,\psi_j)$.

Our main results are given as follows:

\begin{theorem}\label{thm1}
Assume that $\mathrm{VD}$, $\mathrm{RVD}$ and the $\rho$-chain condition hold, and $(\mathcal{E},\mathcal{F})$ is a strongly local regular Dirichlet form satisfying $\mathrm{PHI}(\psi_c)$. Then, the heat kernel $p_t$ exists, and there are constants $C,c_1,c_2>0$ such that for all $x,y\in M$ and $0<t<\psi_c(diam(M,d))$,
\begin{equation}\label{rho-gauss}
p_t(x,y)\asymp\frac{1}{V(x,\psi_c^{-1}(t))}\exp\left\{-c_\bullet\sup\limits_{r>0}\left(\frac{\rho(d(x,y))}{\rho(r)}-\frac{t}{\psi_c(r)}\right)\right\}.
\end{equation}
\end{theorem}

\begin{theorem}\label{thm2}
Assume that $\mathrm{VD}$, $\mathrm{RVD}$ and the $\rho$-chain condition hold, and that there exists a regular strongly local Dirichlet form on $L^2(M,\mu)$ satisfying $\mathrm{PHI}(\psi_c)$. Let $J$ be a jump kernel on $(M,\mu)$ satisfying $\mathrm{J}(\psi_j)$. If (\ref{scl-int}) holds and $\beta_1^{(r_*)}(\Phi\circ\rho^{-1})>1$ for some $r_*\in(0,\infty]$, then $J$ defines a regular jump-type Dirichlet form that admits a heat kernel $p_t$, and there exist $C,c_1,c_2>0$ such that for all $x,y\in M$ and $0<t<\Phi(r_*\wedge diam(M,d))$,
\begin{equation}\label{S+HK}
p_t(x,y)\asymp\frac{1}{V(x,\Phi^{-1}(t))}\wedge\left\{\frac{t}{V(x,d(x,y))\psi_j(d(x,y))}+\frac{1}{V(x,\Phi^{-1}(t))}\exp\left(-c_\bullet\sup\limits_{\sigma>0}\left[\frac{\rho(d(x,y))}{\rho(\sigma)}-\frac{t}{\Phi(\sigma)}\right]\right)\right\}.
\end{equation}
In particular, if $\beta_1(\Phi\circ\rho^{-1})>1$, then (\ref{S+HK}) holds for all $x,y\in M$ and $0<t<diam(M,d)$.
\end{theorem}

Note that on the contrary, (\ref{rho-gauss}) implies $\mathrm{PHI}(\psi_c)$ and (\ref{S+HK}) implies $\mathrm{J}(\psi_j)$ completely in the same way as in \cite{gt} and \cite{bkkl}. Therefore, what we obtain is actually equivalent conditions for these (sharp) heat kernel estimates.

The condition $\beta_1^{(r_*)}(\Phi\circ\rho^{-1})>1$ is satisfied when, for example, $\beta_1^{(r_*)}(\Phi)>\gamma_2$ (see Lemma \ref{beta-uni}). Meanwhile, we do not assume any relation between $\beta_1^{(r_*)}(\Phi)$ (or $\beta_1^{(r_*)}(\psi_j)$, or $\beta_1^{(r_*)}(\psi_c)$) and $\gamma_2$ in general. Note also that when the chain condition holds, that is, when $\rho=\rho_1$, then the estimate (\ref{S+HK}) is equivalent to SHK. In this way, (\ref{S+HK}) can be considered as a generalization of SHK. Further, in our proof, (\ref{S+HK}) is deduced from G$^+$HK under a metric $d_\varepsilon$ and a metric transform from $d_\varepsilon$ to $d$. Additionally, (\ref{S+HK}) implies $\mathrm{PHI}(\Phi)$ in the same way as in \cite[Corollary 2.5]{lm}. If $\Phi\asymp\psi_j$ (which happens if $\psi_j$ is away from $\psi_c$, like $\beta_2(\psi_j\circ\psi_c^{-1})<1$), then all SHK, G$^+$HK and (\ref{S+HK}) are equivalent to the stable-like estimate.

\section{The scale $\rho$ and strongly local Dirichlet forms}

We begin with some properties of a general $\rho$-chain condition:

\begin{lemma}\label{dense-rad}
Assume the $\rho$-chain condition holds on $(M,d)$ with some scale function $\rho$. For any $x\in M$, the set
$$\mathcal{R}_x:=\{r>0:\mbox{ there exists }\xi\in M\mbox{ such that }d(x,\xi)=r\}$$
is dense in $(0,\sup_{\xi\in M}d(x,\xi))$.
\end{lemma}

\begin{proof}
Suppose there exists $0<r_0<\sup_{\xi\in M}d(x,\xi)$ and $0<\delta<r_0\wedge(\sup_{\xi\in M}d(x,\xi)-r_0)$ such that for any $r\in(r_0-\delta,r_0+\delta)$, there is no $\xi\in M$ satisfying $d(x,\xi)=r$. Therefore, $B(x,r_0+\delta)\setminus B(x,r_0-\delta)=\emptyset.$ But from the selection of $r_0$ and $\delta$ we see there exists $y\in M\setminus B(x,r_0+\delta)$, and for $\varepsilon<\delta$, there is no $\varepsilon$-chain from $x$ to $y$. That is, $d_\varepsilon(x,y)=\infty$. However, by (\ref{scale-comp}) it is clear that $d_\varepsilon(x,y)<\infty$ for any $x,y\in M$, which is a contradiction. Thus $\mathcal{R}_x$ is dense in $(0,\sup_{\xi\in M}d(x,\xi))$.
\end{proof}

Given any $\varepsilon>0$, we denote
$$F_{\rho,\varepsilon}(r):=\begin{cases}r,&\mbox{if }r\in(0,\varepsilon];\\
\rho^{-1}\left(\frac{r}{\varepsilon}\rho(\varepsilon)\right),&\mbox{if }r\in[\varepsilon,\infty].
\end{cases}$$
This transform will be frequently used in Section 3.

The following two lemmas give some interesting properties of $F_{\rho,\varepsilon}$:

\begin{lemma}\label{beta-uni}
Let $\psi$ be an arbitrary scale function, and fix $r_*\in(0,\infty]$. Then, $\psi\circ F_{\rho,\varepsilon}$ is a scale function with
$$\beta_1^{(F_{\rho,\varepsilon}^{-1}(r_*))}(\psi\circ F_{\rho,\varepsilon})\ge\beta_1^{(r_*)}(\psi\circ\rho^{-1})\ge\beta_1^{(r_*)}(\psi)/\gamma_2,\quad\beta_2(\psi\circ F_{\rho,\varepsilon})\le\beta_2(\psi).$$
\end{lemma}

\begin{proof}
Let $C=C(\psi)\vee C(d,\rho)$, $\beta'_1=\beta_1^{(r_*)}(\psi)$ and $\beta'_2=\beta_2(\psi)$. For all $0<r_1\le r_2\le\varepsilon$, we have
\begin{equation}\label{scale0}
C^{-1}\left(\frac{r_2}{r_1}\right)^{\beta'_1}\le\frac{\psi\circ F_{\rho,\varepsilon}(r_2)}{\psi\circ F_{\rho,\varepsilon}(r_1)}=\frac{\psi(r_2)}{\psi(r_1)}\le C\left(\frac{r_2}{r_1}\right)^{\beta'_2};
\end{equation}
for all $\varepsilon\le r_1\le r_2<F_{\rho,\varepsilon}(r_*)$, we have
\begin{align*}
(C')^{-1}\left(\frac{r_2}{r_1}\right)^{\beta'_1/\gamma_2}=&C^{-1}\left(C^{-1}\left(\frac{\frac{r_2}{\varepsilon}\rho(\varepsilon)}{\frac{r_1}{\varepsilon}\rho(\varepsilon)}\right)^{1/\gamma_2}\right)^{\beta'_1}\le C^{-1}\left(\frac{\rho^{-1}\left(\frac{r_2}{\varepsilon}\rho(\varepsilon)\right)}{\rho^{-1}\left(\frac{r_1}{\varepsilon}\rho(\varepsilon)\right)}\right)^{\beta'_1}\le\frac{\psi\circ F_{\rho,\varepsilon}(r_2)}{\psi\circ F_{\rho,\varepsilon}(r_1)}=\\
=&\frac{\psi\circ\rho^{-1}\left(\frac{r_2}{\varepsilon}\rho(\varepsilon)\right)}{\psi\circ\rho^{-1}\left(\frac{r_1}{\varepsilon}\rho(\varepsilon)\right)}\le C\left(\frac{\rho^{-1}\left(\frac{r_2}{\varepsilon}\rho(\varepsilon)\right)}{\rho^{-1}\left(\frac{r_1}{\varepsilon}\rho(\varepsilon)\right)}\right)^{\beta'_2}\le C\left(C\left(\frac{\frac{r_2}{\varepsilon}\rho(\varepsilon)}{\frac{r_1}{\varepsilon}\rho(\varepsilon)}\right)^{1/\gamma_1}\right)^{\beta'_2}=C'\left(\frac{r_2}{r_1}\right)^{\beta'_2/\gamma_1};
\end{align*}
while for all $0<r_1<\varepsilon\le r_2<F_{\rho,\varepsilon}(r_*)$,
\begin{align*}
(CC')^{-1}\left(\frac{r_2}{r_1}\right)^{\beta'_1/\gamma_2}\le&(C')^{-1}\left(\frac{r_2}{\varepsilon}\right)^{\beta'_1/\gamma_2}\cdot C^{-1}\left(\frac{\varepsilon}{r_1}\right)^{\beta'_1}\le\frac{\psi\circ F_{\rho,\varepsilon}(r_2)}{\psi\circ F_{\rho,\varepsilon}(r_1)}=\\
=&\frac{\psi\circ\rho^{-1}\left(\frac{r_2}{\varepsilon}\rho(\varepsilon)\right)}{\psi\circ\rho^{-1}(\rho(\varepsilon))}\cdot\frac{\psi(\varepsilon)}{\psi(r_1)}\le C'\left(\frac{r_2}{\varepsilon}\right)^{\beta'_2/\gamma_1}\cdot C^{-1}\left(\frac{\varepsilon}{r_1}\right)^{\beta'_2}\le CC'\left(\frac{r_2}{r_1}\right)^{\beta'_2}.
\end{align*}
Therefore, $\psi\circ F_{\rho,\varepsilon}$ is a scale function with $\beta_1(\psi\circ F_{\rho,\varepsilon})\ge\beta'_1/\gamma_2$ and $\beta_2(\psi\circ F_{\rho,\varepsilon})\le\beta'_2$.
\end{proof}

In particular, $\beta_1^{(r_*)}(F_{\rho,\varepsilon})\ge 1/\gamma_2$ and $\beta_2(F_{\rho,\varepsilon})\le 1$ (since $F_{\rho,\varepsilon}=\rho_1\circ F_{\rho,\varepsilon}$ and clearly $\beta_1(\rho_1)=\beta_2(\rho_1)=1$, where $\rho_1:r\mapsto r$).

\begin{lemma}\label{rho-lem}
If the $\rho$-chain condition holds, then there exists a constant $C\ge 1$ (only depending on the scale function $\rho$) such that for all $x,y\in M$ and $\varepsilon>0$, we have $C^{-1}d(x,y)\le F_{\rho,\varepsilon}(d_\varepsilon(x,y))\le Cd(x,y)$.
\end{lemma}

\begin{proof}
When $\varepsilon>d(x,y)$, the conclusion is trivial, since $F_{\rho,\varepsilon}(d_\varepsilon(x,y))=F_{\rho,\varepsilon}(d(x,y))=d(x,y)$.

When $\varepsilon\le d(x,y)$, by (\ref{rho-chain}) we have
$$\rho^{-1}\left(C^{-1}\frac{\rho(d(x,y))}{\rho(\varepsilon)}\rho(\varepsilon)\right)\le F_{\rho,\varepsilon}(d_\varepsilon(x,y))=\rho^{-1}\left(\frac{d_\varepsilon(x,y)}{\varepsilon}\rho(\varepsilon)\right)\le\rho^{-1}\left(C\frac{\rho(d(x,y))}{\rho(\varepsilon)}\rho(\varepsilon)\right).$$
That is, for all $x,y\in M$ and all $\varepsilon\in(0,d(x,y))$,
$$(C')^{-1}d(x,y)\le\rho^{-1}(C^{-1}\rho(d(x,y)))\le F_{\rho,\varepsilon}(d_\varepsilon(x,y))\le\rho^{-1}(C\rho(d(x,y)))\le C'd(x,y).$$
Therefore, the desired conclusion holds for both cases.
\end{proof}

\begin{corollary}\label{rho-cor}
If the $\rho$-chain condition holds, then there exists a constant $C\ge 1$ (only depending on the scale function $\rho$) such that for all $x,y\in M$ and $\varepsilon>0$, we have
$$C^{-1}d_\varepsilon(x,y)\le F_{\rho,\varepsilon}^{-1}(d(x,y))\le Cd_\varepsilon(x,y).$$
\end{corollary}

\begin{proof}
By (\ref{scale''}) with the scale function $\rho$.
\end{proof}

Now we focus on a metric measure space equipped with a strongly local regular Dirichlet form satisfying $\mathrm{PHI}(\psi_c)$. We remark first that if $\mathrm{PHI}(\psi_c)$ holds for some diffusion, then $\beta_1:=\beta_1(\psi_c)>1$ always holds. See \cite[the proof of Theorem 2.2]{lm}.

\begin{proposition}\label{p-phi-hk}
Let $(M,d,\mu)$ be a metric measure space and $(\mathcal{E}^{(c)},\mathcal{F}^{(c)})$ be a regular strongly local Dirichlet form satisfying $\mathrm{PHI}(\psi_c)$, then there exist constants $C\ge 1$ and $\kappa_1,\kappa_2>0$ such that
\begin{equation}\label{general-loc-hk}
\frac{C^{-1}}{V(x,\psi_c^{-1}(t))}\exp\left\{-C\frac{d_{\varepsilon(\kappa_1t,x,y)}(x,y)}{\varepsilon(\kappa_1t,x,y)}\right\}\le p_t(x,y)\le\frac{C}{V(x,\psi_c^{-1}(t))}\exp\left\{-C^{-1}\frac{d_{\varepsilon(\kappa_2t,x,y)}(x,y)}{\varepsilon(\kappa_2t,x,y)}\right\},
\end{equation}
for all $x,y\in M$ and $0<t<\psi_c(\overline{R})$, where
$$\varepsilon(t,x,y)=\sup\left\{\varepsilon>0: \frac{\psi_c(\varepsilon)}{\varepsilon}d_\varepsilon(x,y)\le t\right\}.$$
\end{proposition}

\begin{proof}
By Proposition \ref{PHI}, $\mathrm{PHI}(\psi_c)\Rightarrow\mathrm{PI}(\psi_c)+\mathrm{cap}(\psi_c)$, and thus by \cite[Theorem 1.6]{mu},
\begin{equation}\label{scale-comp}
\left(\frac{d_\varepsilon(x,y)}{\varepsilon}\right)^2\le C\frac{\psi_c(d(x,y))}{\psi_c(\varepsilon)},\mbox{ for all }\varepsilon<d(x,y).
\end{equation}
and consequently $\lim\limits_{\varepsilon\to 0}\frac{\psi_c(\varepsilon)}{\varepsilon}d_\varepsilon(x,y)=0\mbox{, for any }x,y\in M.$ Therefore, by \cite[Theorem 6.5]{gt}, $\mathrm{VD}+\mathrm{EHI}+\mathrm{E}(\psi_c)\Rightarrow$(\ref{general-loc-hk}). Thus the proposition follows.
\end{proof}

Note by \cite[Lemma 6.4]{gt} that $\varepsilon=\varepsilon(t,x,y)$ satisfies
\begin{equation}\label{eps-value}
\frac{\psi_c(\varepsilon)}{\varepsilon}d_\varepsilon(x,y)=t.
\end{equation}

The inequality (\ref{scale-comp}) indicates some important properties of the scale $\rho$ if the $\rho$-chain condition holds:

\begin{lemma}\label{density}
If the $\rho$-chain condition and (\ref{scale-comp}) holds, then there exists $C_5>0$ such that
\begin{equation}\label{scale-comp+}
\frac{\rho(r_2)}{\rho(r_1)}\le C_5\sqrt{\frac{\psi_c(r_2)}{\psi_c(r_1)}}\quad\mbox{for all }0<r_1\le r_2<\overline{R}.
\end{equation}
\end{lemma}

\begin{proof}
Fix arbitrary $0<r_1\le r_2<\overline{R}$. Then there exists $x_0,y_0\in M$ such that $d(x_0,y_0)>r_2$. In particular,
$$\sup_{\xi\in M}d(x_0,\xi)\ge d(x_0,y_0)>r_2.$$
Thus by Lemma \ref{dense-rad}, for any $\delta>0$, there exists $y\in M$ such that $r_2-\delta<d(x_0,y)<r_2+\delta$, and by (\ref{scale-comp}),
$$\frac{\rho(d(x_0,y))}{\rho(r_1)}\le C\frac{d_\varepsilon(x_0,y)}{r_1}\le C\sqrt{C\frac{\psi_c(d(x_0,y))}{\psi_c(r_1)}}.$$
By the continuity of $\rho$ and $\psi_c$, as $d(x_0,y)\to r_2$, we have
$$\frac{\rho(r_2)}{\rho(r_1)}\le C\sqrt{C\frac{\psi_c(r_2)}{\psi_c(r_1)}}.$$
The proof is complete by letting $C_5=C\sqrt{C}$.
\end{proof}

Note that the value $\rho(r)$, with any $r\ge\overline{R}$, is of no use. Therefore, after proper modification on $[\overline{R},\infty)$ (if not empty), we assume that (\ref{scale-comp+}) holds for all $0<r_1\le r_2<\infty$. As a consequence,

\begin{corollary}\label{beta-uni'}
If $\mathrm{PHI}(\psi_c)$ holds for some regular strongly local Dirichlet form, then $\beta_1(\psi_c\circ F_{\rho,\varepsilon})\ge\beta_1\wedge 2.$
\end{corollary}

\begin{proof}
By Lemma \ref{density}, for all $0<s_1\le s_2$, letting $\sigma=\rho^{-1}(s_i)$, we have
$$\frac{s_2}{s_1}=\frac{\rho(\sigma_2)}{\rho(\sigma_1)}\le C_5\sqrt{\frac{\psi_c(\sigma_2)}{\psi_c(\sigma_1)}}=C_5\sqrt{\frac{\psi_c\circ\rho^{-1}(s_2)}{\psi_c\circ\rho^{-1}(s_1)}}.$$
Therefore, for all $\varepsilon\le r_1\le r_2$, we have
$$\frac{\psi_c\circ F_{\rho,\varepsilon}(r_2)}{\psi_c\circ F_{\rho,\varepsilon}(r_1)}=\frac{\psi_c\circ\rho^{-1}\left(\frac{r_2}{\varepsilon}\rho(\varepsilon)\right)}{\psi_c\circ\rho^{-1}\left(\frac{r_1}{\varepsilon}\rho(\varepsilon)\right)}\ge C_5^{-2}\left(\frac{\frac{r_2}{\varepsilon}\rho(\varepsilon)}{\frac{r_1}{\varepsilon}\rho(\varepsilon)}\right)^2=C_5^{-2}\left(\frac{r_2}{r_1}\right)^2;$$
(\ref{scale0}) holds for $0<r_1\le r_2\le\varepsilon$; while for all $0<r_1<\varepsilon\le r_2$,
$$\frac{\psi_c\circ F_{\rho,\varepsilon}(r_2)}{\psi_c\circ F_{\rho,\varepsilon}(r_1)}=\frac{\psi_c\circ\rho^{-1}\left(\frac{r_2}{\varepsilon}\rho(\varepsilon)\right)}{\psi_c\circ\rho^{-1}(\rho(\varepsilon))}\cdot\frac{\psi_c(\varepsilon)}{\psi_c(r_1)}\ge C_5^{-2}\left(\frac{r_2}{\varepsilon}\right)^2\cdot C^{-1}\left(\frac{\varepsilon}{r_1}\right)^{\beta_1}\ge(CC_5^2)^{-1}\left(\frac{r_2}{r_1}\right)^{\beta_1\wedge 2}.$$
Therefore, $\beta_1(\psi_c\circ F_{\rho,\varepsilon})\ge\beta_1\wedge 2$.
\end{proof}

\begin{proof}[Proof of Theorem \ref{thm1}]
By Lemma \ref{3.8} below (clearly independent of Section 2), it suffices to prove that there exist constants $C,c_1,c_2,\kappa_1,\kappa_2>0$ such that
$$\frac{C^{-1}}{V(x,\psi_c^{-1}(t))}\exp\left\{-c_1\sup\limits_{r>0}\left(\frac{\rho(d(x,y))}{\rho(r)}-\frac{\kappa_1t}{\psi_c(r)}\right)\right\}\le p_t(x,y)\le\frac{C}{V(x,\psi_c^{-1}(t))}\exp\left\{-c_2\sup\limits_{r>0}\left(\frac{\rho(d(x,y))}{\rho(r)}-\frac{\kappa_2t}{\psi_c(r)}\right)\right\}.$$

Let $f_t(r)=f_t(r;x,y):=\frac{\rho(d(x,y))}{\rho(r)}-\frac{t}{\psi_c(r)}$. By (\ref{eps-value}) we have
$$\sup\limits_{r>0}f_t(r)=\sup\limits_{r>0}\left(\frac{\rho(d(x,y))}{\rho(r)}-\frac{t}{\psi_c(r)}\right)\le\frac{d_\varepsilon(x,y)}{\varepsilon}\sup\limits_{r>0}\left(C\frac{\rho(\varepsilon)}{\rho(r)}-\frac{\psi_c(\varepsilon)}{\psi_c(r)}\right).$$
where $C=C(d,\rho)$. Note that if $r\ge\varepsilon$,
$$C\frac{\rho(\varepsilon)}{\rho(r)}-\frac{\psi_c(\varepsilon)}{\psi_c(r)}\le C\cdot 1-0=C,$$
while if $0<r<\varepsilon$, by (\ref{scale-comp+}) we have
$$C\frac{\rho(\varepsilon)}{\rho(r)}-\frac{\psi_c(\varepsilon)}{\psi_c(r)}\le CC_5\sqrt{\frac{\psi_c(\varepsilon)}{\psi_c(r)}}-\frac{\psi_c(\varepsilon)}{\psi_c(r)}\le\frac{1}{4}(CC_5)^2.$$

On the contrary,
$$\sup\limits_{r>0}f_t(r)\ge\sup\limits_{r\ge\varepsilon}\left(\frac{\rho(d(x,y))}{\rho(r)}-\frac{t}{\psi_c(r)}\right)\ge\frac{d_\varepsilon(x,y)}{\varepsilon}\sup\limits_{r\ge\varepsilon}\left(C^{-1}\frac{\rho(\varepsilon)}{\rho(r)}-\frac{\psi_c(\varepsilon)}{\psi_c(r)}\right).$$
Again by (\ref{scale-comp+}) we have
$$\sup\limits_{r\ge\varepsilon}\left(C^{-1}\frac{\rho(\varepsilon)}{\rho(r)}-\frac{\psi_c(\varepsilon)}{\psi_c(r)}\right)\ge\sup\limits_{r\ge\varepsilon}\left(C^{-1}C_5^{-1}\sqrt{\frac{\psi_c(\varepsilon)}{\psi_c(r)}}-\frac{\psi_c(\varepsilon)}{\psi_c(r)}\right)=\frac{1}{4}(CC_5)^{-2}.$$

Therefore,
\begin{equation}\label{opt-n_eps}
\frac{1}{4}(CC_5)^{-2}\cdot\frac{d_\varepsilon(x,y)}{\varepsilon}\le\sup\limits_{r>0}f_t(r)\le \left(C\vee\frac{1}{4}(CC_5)^2\right)\frac{d_\varepsilon(x,y)}{\varepsilon}.
\end{equation}
Note also that by Proposition \ref{p-phi-hk}, the heat kernel estimate (\ref{general-loc-hk}) holds. Therefore,
\begin{align*}
\frac{C^{-1}}{V(x,\psi_c^{-1}(t))}\exp\left\{-4C(CC_5)^2\sup\limits_{r>0}\left(\frac{\rho(d(x,y))}{\rho(r)}-\frac{\kappa_1t}{\psi_c(r)}\right)\right\}\le\frac{C^{-1}}{V(x,\psi_c^{-1}(t))}\exp\left\{-C\frac{d_{\varepsilon(\kappa_1t,x,y)}(x,y)}{\varepsilon(\kappa_1t,x,y)}\right\}\le p_t(x,y)\\
\le\frac{C}{V(x,\psi_c^{-1}(t))}\exp\left\{-C^{-1}\frac{d_{\varepsilon(\kappa_2t,x,y)}(x,y)}{\varepsilon(\kappa_2t,x,y)}\right\}\le\frac{C}{V(x,\psi_c^{-1}(t))}\exp\left\{\left(C^{-2}\wedge\frac{4C^{-1}}{(CC_5)^2}\right)\sup\limits_{r>0}\left(\frac{\rho(d(x,y))}{\rho(r)}-\frac{\kappa_2t}{\psi_c(r)}\right)\right\}.
\end{align*}
The theorem follows by letting
$$c_1=4C(CC_5)^2,\quad c_2=C^{-2}\wedge\frac{4C^{-1}}{(CC_5)^2},$$
and $\kappa_1,\kappa_2$ come from Proposition \ref{p-phi-hk}.
\end{proof}

From this theorem we can prove $\mathrm{NLE}(\psi_c)$ in a strong sense (that is, $\eta=1$):

\begin{corollary}\label{diff-nle}
Assume that the $\rho$-chain condition holds. For a diffusion satisfying $\mathrm{PHI}(\psi_c)$, there exists $C_6>0$ such that for all $x\in M$ and $0<t<\psi_c(\overline{R})$,
\begin{equation}\label{nle}
p_t(x,y)\ge\frac{C_6^{-1}}{V(x,\psi_c^{-1}(t))}\mbox{ on }B(x,\psi_c^{-1}(t)).
\end{equation}
\end{corollary}

\begin{proof}
Note that under condition PHI$(\psi_c)$, (\ref{scale-comp+}) holds. Therefore,
\begin{align*}
\sup\limits_{r>0}\left(\frac{\rho(d(x,y))}{\rho(r)}-\frac{\kappa_1t}{\psi_c(r)}\right)=&\sup\limits_{0<r<d(x,y)}\left(\frac{\rho(d(x,y))}{\rho(r)}-\frac{\kappa_1t}{\psi_c(r)}\right)\vee\sup\limits_{r\ge d(x,y)}\left(\frac{\rho(d(x,y))}{\rho(r)}-\frac{\kappa_1t}{\psi_c(r)}\right)\\
\le&\sup\limits_{0<r<d(x,y)}\left(C_5\sqrt{\frac{\psi_c(d(x,y))}{\psi_c(r)}}-\frac{\kappa_1t}{\psi_c(r)}\right)\vee\left(\frac{\rho(d(x,y))}{\rho(d(x,y))}-0\right)=\frac{C_5^2\psi_c(d(x,y))}{4\kappa_1t}\vee 1.
\end{align*}
Thus by (\ref{rho-gauss}), the estimate (\ref{nle}) holds with $C_6=C\exp(c_1C_5^2/(4\kappa_1))$.
\end{proof}

\section{Jump-type Dirichlet forms}

Here are some basic properties of $\Phi$ defined by (\ref{phi}), from which we see easily that $\Phi$ is a scale function:

\begin{lemma}\cite[Corollary 2.5]{lm}
Provided that $\int_0^1\frac{\psi_c(s)ds}{s\psi_j(s)}<\infty$, there exists a constant $C>0$ such that:
\begin{align}
\Phi(r)\le C\psi_j(r),&\quad\mbox{for all }r>0;\label{comp2}\\
\frac{\Phi(R)}{\Phi(r)}\le C\frac{\psi_c(R)}{\psi_c(r)},&\quad\mbox{for all }0<r\le R;\label{comp3}\\
\Phi(r)\le C\psi_c(r),&\quad\mbox{for all }r\ge 1;\notag\\
\psi_c(r)\le C\Phi(r),&\quad\mbox{for all }0<r\le 1.\notag
\end{align}
\end{lemma}

Note by \cite[Lemma 3.5]{lm} that the condition $\int_0^1\frac{\psi_c(s)ds}{s\psi_j(s)}<\infty$ here is equivalent to (\ref{scl-int}).

The following proposition, which comes from \cite{bkkl}, is prerequisite for our work:

\begin{proposition}\label{found}
Assume that $(M,d)$ is unbounded, satisfying $\mathrm{VD}$ and $\mathrm{RVD}$. Assume also that there exists a regular strongly local Dirichlet form satisfying $\mathrm{PHI}(\psi_c)$. If $\psi_j$ is a scale function such that $\Phi$ (defined by (\ref{phi})) is a scale function with $\beta_1^{(T)}(\Phi)>1$, where $T\in(0,\infty]$, then $\mathrm{J}(\psi_j)\Leftrightarrow\mathrm{J}(\psi_j)+\mathrm{PI}(\Phi)+\mathrm{E}(\Phi)\Leftrightarrow\mathrm{G^+HK}(\Phi,\psi_j;T).$ Further, if the chain condition holds, then $\mathrm{J}(\psi_j)\Leftrightarrow\mathrm{SHK}(\Phi,\psi_j;T).$
\end{proposition}

\begin{proof}
$\mathrm{J}(\psi_j)\Rightarrow\mathrm{PI}(\Phi)$ by \cite[Theorem 4.6]{bkkl}, $\mathrm{J}(\psi_j)+\mathrm{PI}(\Phi)+\mathrm{E}(\Phi)\Rightarrow\mathrm{G^+HK}(\Phi,\psi_j)$ by \cite[Corollary 2.15 and Proposition 3.12]{bkkl}, and $\mathrm{G^+HK}(\Phi,\psi_j)\Rightarrow\mathrm{J}(\psi_j)$ is trivial. The equivalence under the chain condition comes from \cite[Theorem 2.17]{bkkl}.
\end{proof}

\subsection{The proof to $\mathrm{UHK}(\Phi,\psi_j)$ revisited}

We will prove Theorem \ref{thm2} based on the result above, but the version in Proposition \ref{found} is not sufficient. To strengthen it, we introduce the uniform stability of $\mathrm{UHK}(\Phi,\psi_j)$. Let $I$ be an arbitrary index set. We say

\begin{itemize}
\item real numbers $r_i>1$ uniformly on $i\in I$, if there exists $\varepsilon>0$ such that for each $i$, $r_i\ge 1+\varepsilon$;

\item functions $\phi^{(i)}$ are scale functions uniformly on $i\in I$ (with parameters $\beta_{1,0},\beta_{2,0},C_{\phi,0}$), if $0<\beta_{1,0}\le\beta_{2,0}<\infty$ and $C_{\phi,0}>1$ satisfy $\beta_1(\phi^{(i)})\ge\beta_{1,0},\beta_2(\phi^{(i)})\le\beta_{2,0}$ and $C(\phi^{(i)})\le C_{\phi,0}$;

\item functions $\phi^{(i)}$ are scale functions strongly uniformly on $i\in I$ (with parameters $\beta_{*,0},\beta_{1,0},\beta_{2,0},C_{\phi,0}$), if $0<\beta_{1,0}\le\beta_{2,0}<\infty$, $\beta_{*,0},C_{\phi,0}>1$ satisfy $\beta_1(\phi^{(i)})\ge\beta_{1,0},\beta_2(\phi^{(i)})\le\beta_{2,0}$, $C(\phi^{(i)})\cdot C'(\phi^{(i)})\le C_{\phi,0}$ and $\beta_1^{(r_{*,i})}(\phi^{(i)})\ge\beta_{*,0}$ for some $r_{*,i}\in(0,\infty]$;

\item $(M^{(i)},d^{(i)},\mu^{(i)})$ satisfy $\mathrm{RVD}$ (or respectively $\mathrm{VD}$) uniformly on $i\in I$, if there exist $\alpha_{1,0},C_{RVD,0}\ge 1$ (respectively, $\alpha_{2,0},C_{VD,0}\ge 1$) such that for all $i\in I$, $\mathrm{RVD}$ holds with $\alpha_1(M^{(i)},d^{(i)},\mu^{(i)})\ge\alpha_{1,0}$ and $C_{RVD}(M^{(i)},d^{(i)},\mu^{(i)})\le C_{RVD,0}$ (respectively, $\mathrm{VD}$ with $\alpha_2(M^{(i)},d^{(i)},\mu^{(i)})\le\alpha_{2,0}$ and $C_{VD}(M^{(i)},d^{(i)},\mu^{(i)})\le C_{VD,0}$);

\item jump kernels $J^{(i)}$ satisfy $\mathrm{J}(\psi_j^{(i)})$ uniformly on $i\in I$, if $\psi_j^{(i)}$ are scale functions uniformly and there exists $C_{J,0}\ge 1$ such that for all $i\in I$, $\mathrm{J}(\psi_j^{(i)})$ holds with $C_J(J,\psi_j^{(i)})\le C_{J,0}$;

\item Dirichlet spaces $(M^{(i)},d^{(i)},\mu^{(i)},\mathcal{E}^{(i)},\mathcal{F}^{(i)})$ satisfy $\mathrm{E}(\phi^{(i)})$ (or respectively $\mathrm{DUE}(\phi^{(i)})$) uniformly on $i\in I$, if $\phi^{(i)}$ are scale functions uniformly and there exist $C_{E,0}\ge 1$ and $\kappa_{E,0}>0$ (respectively $C_{DUE,0}\ge 1$) such that for all $i\in I$, $\mathrm{E}(\phi^{(i)})$ holds with $C_E(\mathcal{E}^{(i)},\mathcal{F}^{(i)},\phi^{(i)})\le C_{E,0}$ and $\kappa_E(\mathcal{E}^{(i)},\mathcal{F}^{(i)},\phi^{(i)})\ge \kappa_{E,0}$ (respectively $\mathrm{DUE}(\phi^{(i)})$ with $C_{DUE}(\mathcal{E}^{(i)},\mathcal{F}^{(i)},\phi^{(i)})\le C_{DUE,0}$).
\end{itemize}

Here $\beta_{*,0}>1$ is extremely important for the strong uniformity of scale functions. Note that when $r_{*,i}=\infty$ for all $i\in I$, we can (and always do) take $\beta_{1,0}=\beta_{*,0}>1$, and the (ordinary) uniformity follows.

Now we introduce some notions for simplicity, which are only valid in this subsection. For any $i\in I$, $x_i,y_i\in M^{(i)}$ and $r>0$, let $\psi_j^{-i}:=(\psi_j^{(i)})^{-1}$, $\phi^{-i}:=(\phi^{(i)})^{-1}$, and
$$B_i(x_i,r):=B_{d^{(i)}}(x_i,r),\ V^{(i)}(x_i,r):=\mu^{(i)}(B_i(x_i,r)),\ V^{(i)}(x_i,y_i):=V^{(i)}(x_i,d^{(i)}(x_i,y_i)),\ \psi_j^{(i)}(x_i,y_i):=\psi_j^{(i)}(d^{(i)}(x_i,y_i)).$$
Let $\{P^{(i)}_t\}$ be the heat semigroup corresponding to $(M^{(i)},d^{(i)},\mu^{(i)},\mathcal{E}^{(i)},\mathcal{F}^{(i)})$, and $p^{(i)}_t(x,y)$ be its heat kernel. If $\mathcal{E}^{(i)}$ is a pure jump-type form with a jump kernel $J^{(i)}$, let
$$\mathcal{E}^{(i,\sigma)}(u):=\int_{M^{(i)}}\int_{M^{(i)}}(u(x)-u(y))^2J^{(i)}(x,y)d\mu^{(i)}(y)d\mu^{(i)}(x)$$
be its $\sigma$-truncated form. If $(\mathcal{E}^{(i,\sigma)},\mathcal{F}^{(i)})$ is also a regular Dirichlet form (which is always true provided that $\mathrm{J}_\le(\psi)$ holds with an arbitrary scale function $\psi$, see \cite[Proposition 4.2 and Remark 4.3]{ghl14}), let $\{P^{(i,\sigma)}_t\}$ be its heat semigroup and $p^{(i,\sigma)}_t(x,y)$ be its heat kernel. Further, given any open subset $\Omega\subset M^{(i)}$, let $P^{(i,\Omega)}_t:=(P^{(i)})^\Omega_t$ and $P^{(i,\sigma,\Omega)}_t:=(P^{(i,\sigma)})^\Omega_t$.

\begin{proposition}\label{found+}
Let $\{(M^{(i)},d^{(i)},\mu^{(i)})\}_{i\in I}$ be a family of metric measure spaces satisfying $\mathrm{VD}$ uniformly, and $J^{(i)}$ be a jump kernel on $M^{(i)}$ satisfying $\mathrm{J}(\psi_j^{(i)})$ uniformly on $i$, where $\psi_j^{(i)}$ are scale functions uniformly with parameters $\beta_{1,0},\beta_{2,0},C_{\psi,0}$. Assume that each $J^{(i)}$ defines a regular Dirichlet form $(\mathcal{E}^{(i)},\mathcal{F}^{(i)})$ satisfying $\mathrm{E}(\phi^{(i)})$ and $\mathrm{DUE}(\phi^{(i)})$ uniformly, where $\phi^{(i)}$ are strongly uniform scale functions with parameters $\beta_{*,0},\beta_{1,1},\beta_{2,1},C_{\phi,0}$. Assume also that $\phi^{(i)}(r)\le C_{\phi,1}\psi_j^{(i)}(r)$ with a uniform $C_{\phi,1}\ge 1$. Then $\mathrm{UHK}(\phi^{(i)},\psi_j^{(i)};\phi^{(i)}(r_{*,i}))$ holds uniformly, i.e., there exist constants $C,c>0$ such that for all $i\in I$, $0<t<\phi^{(i)}(r_{*,i}\wedge diam(M^{(i)},d^{(i)}))$ and $x_i,y_i\in M^{(i)}$,
\begin{equation}\label{uhk}
p_t^{(i)}(x_i,y_i)\le\frac{Ct}{V^{(i)}(x_i,y_i)\psi_j^{(i)}(x_i,y_i)}+\frac{C}{V^{(i)}(x_i,\phi^{-i}(t))}\exp\left(-c\sup\limits_{\sigma>0}\left[\frac{d^{(i)}(x_i,y_i)}{\sigma}-\frac{t}{\phi^{(i)}(\sigma)}\right]\right).
\end{equation}
\end{proposition}

Note that (RVD) is not assumed, but is clearly contained in the condition $\mathrm{E}(\phi^{(i)})$ (similar to \cite[Lemma 6.9]{hl}) and is thus automatically uniform on $i\in I$ under the same conditions.

The proof is almost contained in \cite[Subsections 3.1-3.3]{bkkl}, but it is necessary to add more details indicating how the constants (and intermediate parameters) are determined by the given conditions and how boundedness constrains the interval of $t$ when the desired conclusions hold. For this purpose we list the following lemmas (which can be considered as an appendix). The proof of Proposition \ref{found+} will be given after these lemmas. Denote $\overline{R}_i:=diam(M^{(i)},d^{(i)}).$ Let $C_*(\alpha,\gamma,\delta)=\sup\limits_{s\ge 0}(1+2s)^\alpha\exp(-\gamma s^\delta)$ and $c_*(\alpha,\gamma)=\sup\limits_{n\ge 0}(2n-2^{n\alpha}\gamma)_+$ for any $\alpha,\gamma,\delta>0$. Clearly $C_*(\alpha,\gamma,\delta)\le C_*(\alpha',\gamma,\delta)$ if $\alpha<\alpha'$. For simplicity, we assume that $C_{\phi,1}=1$ (otherwise we take $C_{\phi,1}\psi_j^{(i)}$ instead of $\psi_j^{(i)}$ and modify $C_{J,0}$ accordingly). For every step, detailed computation is omitted if it repeats the cited context exactly, while we only mark the difference if there are only slight modifications made on the cited context.

\begin{lemma}\cite[Lemma 3.2 covered]{bkkl}\label{3.2}
Under the same assumptions as in Proposition \ref{found+}, there exist constants $c_6,c_7,c_8,c_9>0$ such that for all $i\in I$, $x_i,y_i\in M^{(i)}$, $\sigma>0$ and $0<t<2\phi^{(i)}(\overline{R}_i)$,
$$p^{(i)}_t(x_i,y_i)\le\frac{c_6}{V^{(i)}(x_i,\phi^{-i}(t))}\left(1+\frac{d^{(i)}(x_i,y_i)}{\phi^{-i}(t)}\right)^{\alpha_{2,0}}\exp\left\{\frac{c_7t}{\phi^{(i)}(\sigma)}-c_8\frac{d^{(i)}(x_i,y_i)}{\sigma}\right\}+\frac{c_9}{V^{(i)}(x_i,\sigma)}\left(\frac{t}{\phi^{(i)}(\sigma)}+1\right)^{\alpha_{2,0}+1}.$$
\end{lemma}
\begin{proof}
By \cite[Lemma 2.1]{ckw0}, for every $i\in I,x_i\in M^{(i)}$ and $r>0$,
$$\int_{B_i(x_i,r)^c}J^{(i)}(x_i,y_i)d\mu^{(i)}(y_i)\le C_{J,0}\sum\limits_{k=0}^\infty\frac{V^{(i)}(x_i,2^{k+1}r)}{V^{(i)}(x_i,2^kr)\psi_j^{(i)}(2^kr)}\le\frac{C_{VD,0}C_{J,0}C_{\psi,0}2^{\alpha_{2,0}}}{1-2^{-\beta_{1,0}}}\frac{1}{\psi_j^{(i)}(r)}=:\frac{c'_5}{\psi_j^{(i)}(r)}.$$
By \cite[Proposition 4.6]{ghl14}, for every open set $\Omega_i\subset M^{(i)}$ and every $\sigma>0$,
\begin{equation}\label{trunc}
\left|P^{(i,\Omega_i)}_tf_i-P^{(i,\sigma,\Omega_i)}_tf_i\right|\le 2t\|f_i\|_{L^\infty(M^{(i)},\mu^{(i)})}\sup\limits_{x_i\in M^{(i)}}\int_{B_i(x_i,\sigma)^c}J^{(i)}(x_i,y_i)d\mu^{(i)}(y_i)\le\frac{2c'_5t}{\psi_j^{(i)}(\sigma)}\|f_i\|_{L^\infty(M^{(i)},\mu^{(i)})}.
\end{equation}

By \cite[Lemma 10.6]{bghh}, given an arbitrary ball $B_i=B_i(x_i,r)$, if $r<\kappa_{E,0}\overline{R}_i$, then for all $t\le(2C_{E,0})^{-1}\phi^{(i)}(r)$ and $\mu^{(i)}$-a.e. $z_i\in\frac{1}{4}B_i$,
\begin{equation}\label{serv}
1-P^{(i,B_i)}_{t}1_{B_i}(z_i)\le 1-\frac{\left(G^{(i,B_i)}1_{B_i}(z_i)-t\right)_+}{\|G^{(i,B_i)}1_{B_i}\|_{L^\infty(M^{(i)},\mu^{(i)})}}\le 1-\frac{\left(C_{E,0}^{-1}\phi^{(i)}(r)-t\right)_+}{C_{E,0}\phi^{(i)}(r)}\le 1-\frac{1}{2C_{E,0}^2}.
\end{equation}
If $r\ge 2\overline{R}_i$, then $1-P_t^{(i,B_i)}1_{B_i}=0$ since $P_t^{(i,B_i)}=P_t^{(i)}$ is conservative by Proposition \ref{PHI}. Thus (\ref{serv}) also holds. Combining (\ref{serv}) with (\ref{trunc}), we have on $\frac{1}{4}B_i$ (where $r<\kappa_{E,0}\overline{R}_i$ or $r\ge 2\overline{R}_i$) with $t\le(2C_{E,0})^{-1}\phi^{(i)}(r)$ that
$$1-P^{(i,\sigma,B_i)}_{t}1_{B_i}\le 1-\frac{1}{2C_{E,0}^2}+2c'_5\frac{t}{\psi_j^{(i)}(\sigma)}\le 1-\frac{1}{2C_{E,0}^2}+2c'_5\frac{t}{\phi^{(i)}(\sigma)}\le 1-\frac{1}{2C_{E,0}^2}+2c'_5\frac{t}{\phi^{(i)}(r\wedge\sigma)}.$$
The same hols for $t>(2C_{E,0})^{-1}\phi^{(i)}(r)$, since
$$1-P^{(i,\sigma,B_i)}_{t}1_{B_i}\le 1<1-\frac{1}{2C_{E,0}^2}+\frac{2c'_5}{2C_{E,0}}\le 1-\frac{1}{2C_{E,0}^2}+2c'_5\frac{(2C_{E,0})^{-1}\phi^{(i)}(r)}{\phi^{(i)}(r\wedge\sigma)}<1-\frac{1}{2C_{E,0}^2}+2c'_5\frac{t}{\phi^{(i)}(r\wedge\sigma)}.$$
If $\kappa_{E,0}\overline{R}_i\le r<2\overline{R}_i$, then applying this result on $B':=B_i(z,\frac{\kappa_{E,0}}{2}\overline{R}_i)$, we have for all $t>0$ and $z\in\frac{1}{4}B_i$ that
$$1-P^{(i,\sigma,B_i)}_{t}1_{B_i}(z)\le 1-P^{(i,\sigma,B')}_{t}1_{B'}(z)\le 1-\frac{1}{2C_{E,0}^2}+2c'_5\frac{t}{\phi^{(i)}(\frac{\kappa_{E,0}}{2}\overline{R}_i\wedge\sigma)}\le 1-\frac{1}{2C_{E,0}^2}+2c'_5\left(\frac{4}{\kappa_{E,0}}\right)^{\beta_{1,0}}C_{\phi,0}\frac{t}{\phi^{(i)}(r\wedge\sigma)}.$$
The last inequality then holds on $\frac{1}{4}B_i$ for all $t,\sigma,r>0$. Thus by \cite[Corollary 4.22]{ckw0}, on $\frac{1}{4}B_i$, for all $t>0$,
\begin{equation}\label{trunc-tail}
P^{(i,\sigma)}_t1_{B_i^c}\le 1-P^{(i,\sigma,B_i)}_t1_{B_i}\le c'_6\exp\left(c'_7\frac{t}{\phi^{(i)}(\sigma)}-c'_8\frac{r}{\sigma}\right),
\end{equation}
where (by taking $\varepsilon=(8C_{E,0}^2)^{-1}$ in \cite[Lemma 4.20]{ckw0} and $\varepsilon'=\varepsilon/2$ in \cite[Lemma 4.21]{ckw0})
$$c'_6=\left(\frac{16C_{E,0}^2-1}{16C_{E,0}^2-2}\right)^2,\quad c'_7=16c'_5C_{\phi,0}C_{E,0}^2\left(\frac{4}{\kappa_{E,0}}\right)^{\beta_{1,0}}\log(8C_{E,0}^2),\quad c'_8=\frac{1}{2}\log\frac{16C_{E,0}^2-1}{16C_{E,0}^2-2}.$$

Further, by \cite[Proposition 3.1]{hl}, for all $x_i,y_i\in M^{(i)}$, $\sigma>0$ and $0<t<\phi^{(i)}(\overline{R}_i)$,
$$P^{(i,\sigma)}_t(x_i,y_i)\le\exp\left(\frac{2c'_5t}{\psi_j^{(i)}(\sigma)}\right)p_t^{(i)}(x_i,y_i)\le\exp\left(\frac{2c'_5t}{\psi_j^{(i)}(\sigma)}\right)\frac{C_{DUE,0}}{V^{(i)}(x_i,\phi^{-i}(t))}.$$
Therefore, by the semigroup property, we obtain for all $0<t<2\phi^{(i)}(\overline{R}_i)$ that
\begin{align}
P^{(i,\sigma)}_t(x_i,y_i)\le%&\left(\int_{B_i(x_i,d^{(i)}(x_i,y_i)/2)^c}+\int_{B_i(y_i,d^{(i)}(x_i,y_i)/2)^c}\right)(p^{(i)})_{t/2}^{(\sigma)}(x_i,z_i)(p^{(i)})_{t/2}^{(\sigma)}(z_i,y_i)d\mu^{(i)}(z_i)\notag\\
&P^{(i,\sigma)}_{t/2}1_{B_i(x_i,d^{(i)}(x_i,y_i)/2)^c}(x_i)\sup\limits_{M^{(i)}}(p^{(i)})_{\frac{t}{2}}^{(\sigma)}(\cdot,y_i)+P^{(i,\sigma)}_{t/2}1_{B_i(y_i,d^{(i)}(x_i,y_i)/2)^c}(y_i)\sup\limits_{M^{(i)}}(p^{(i)})_{\frac{t}{2}}^{(\sigma)}(x_i,\cdot)\notag\\
\le&c'_6\exp\left(c'_7\frac{t/2}{\phi^{(i)}(\sigma)}-c'_8\frac{d^{(i)}(x_i,y_i)}{2\sigma}\right)\exp\left(\frac{2c'_5t/2}{\psi_j^{(i)}(\sigma)}\right)\left(\frac{C_{DUE,0}}{V^{(i)}(x_i,\phi^{-i}(\frac{t}{2}))}+\frac{C_{DUE,0}}{V^{(i)}(y_i,\phi^{-i}(\frac{t}{2}))}\right)\notag\\
\le&\frac{2^{\alpha_{2,0}}c'_6C_{VD,0}^2C_{DUE,0}}{V^{(i)}(x_i,\phi^{-i}(t))}\left(1+\frac{d^{(i)}(x_i,y_i)}{\phi^{-i}(t)}\right)^{\alpha_{2,0}}\exp\left\{\left(c'_5+\frac{c'_7}{2}\right)\frac{t}{\phi^{(i)}(\sigma)}-\frac{c'_8}{2}\frac{d^{(i)}(x_i,y_i)}{\sigma}\right\}.\label{trunc-hk}
\end{align}
Further, parallel to \cite[Lemma 3.1]{bkkl} (where $k_i=\lceil\frac{2c'_7s}{c'_8\phi^{(i)}(\sigma)}\rceil$, noting that the definition $\lceil\cdot\rceil$ is slightly different),
\begin{align*}
\int_{M^{(i)}}&\frac{p^{(i,\sigma)}_s(x_i,z_i)}{V^{(i)}(z_i,\sigma)}d\mu^{(i)}(z_i)=\left(\int_{B_i(x_i,k_i\sigma)}+\sum\limits_{k=k_i}^\infty\int_{B_i(x_i,(k+1)\sigma)\setminus B_i(x_i,k\sigma)}\right)\frac{p^{(i,\sigma)}_s(x_i,z_i)}{V^{(i)}(z_i,\sigma)}d\mu^{(i)}(z_i)\\
\le&\frac{C_{VD,0}(k_i+1)^{\alpha_{2,0}}}{V^{(i)}(x_i,\sigma)}+\sum\limits_{k=k_i}^\infty\frac{C_{VD,0}(k+1+1)^{\alpha_{2,0}}}{V^{(i)}(x_i,\sigma)}c'_6\exp\left(\frac{c'_7s}{\phi^{(i)}(\sigma)}-c'_8k\right)\\
\le&\frac{C_{VD,0}(k_i+1)^{\alpha_{2,0}}}{V^{(i)}(x_i,\sigma)}+\frac{C_{VD,0}c'_6}{V^{(i)}(x_i,\sigma)}\sum\limits_{k=k_i}^\infty(2k+1)^{\alpha_{2,0}}\exp\left(-\frac{c'_8}{2}k\right)\\
\le&\frac{C_{VD,0}}{V^{(i)}(x_i,\sigma)}\left(\frac{2c'_7}{c'_8}\frac{s}{\phi^{(i)}(\sigma)}+3\right)^{\alpha_{2,0}}+\frac{C_{VD,0}c'_6}{V^{(i)}(x_i,\sigma)}\frac{C_*(\alpha_{2,0},c'_8/4,1)}{1-e^{-c'_8/4}}\\
\le&C_{VD,0}\left(c'_6\frac{C_*(\alpha_{2,0},\frac{c'_8}{4},1)}{1-e^{-c'_8/4}}\vee\frac{2c'_7}{c'_8}\vee 3\right)\frac{1}{V^{(i)}(x_i,\sigma)}\left(\frac{s}{\phi^{(i)}(\sigma)}+1\right)^{\alpha_{2,0}}=:\frac{c'_9}{V^{(i)}(x_i,\sigma)}\left(\frac{s}{\phi^{(i)}(\sigma)}+1\right)^{\alpha_{2,0}}.
\end{align*}
Inserting this inequality and (\ref{trunc-hk}) into \cite[Formula (6.4)]{hl}, we obtain
\begin{align}
p^{(i)}_t(x_i,y_i)\le&p^{(i,\sigma)}_t(x_i,y_i)+2\int_0^tds\int_{M^{(i)}}p^{(i)}_{t-s}(z_i,y_i)\int_{B_i(z_i,\sigma)^c}p^{(i,\sigma)}_s(x_i,w_i)J^{(i)}(w_i,z_i)d\mu^{(i)}(w_i)d\mu^{(i)}(z_i)\notag\\
\le&p^{(i,\sigma)}_t(x_i,y_i)+\frac{2C_{J,0}}{\psi_j^{(i)}(\sigma)}\int_0^t\int_{M^{(i)}}p^{(i)}_{t-s}(z_i,y_i)d\mu^{(i)}(z_i)\int_{M^{(i)}}p^{(i,\sigma)}_s(x_i,w_i)\frac{1}{V^{(i)}(w_i,\sigma)}d\mu^{(i)}(w_i)ds\notag\\
\le&p^{(i,\sigma)}_t(x_i,y_i)+\frac{2C_{J,0}}{\psi_j^{(i)}(\sigma)}\int_0^{t}c'_9\frac{1}{V^{(i)}(x_i,\sigma)}\left(\frac{s}{\phi^{(i)}(\sigma)}+1\right)^{\alpha_{2,0}}ds\notag\\
\le&\frac{2^{\alpha_{2,0}}c'_6C_{VD,0}^2C_{DUE,0}}{V^{(i)}(x_i,\phi^{-i}(t))}\left(1+\frac{d^{(i)}(x_i,y_i)}{\phi^{-i}(t)}\right)^{\alpha_{2,0}}\exp\left\{\left(c'_5+\frac{c'_7}{2}\right)\frac{t}{\phi^{(i)}(\sigma)}-\frac{c'_8}{2}\frac{d^{(i)}(x_i,y_i)}{\sigma}\right\}+\notag\\
&\quad\quad\quad+\frac{2C_{J,0}c'_9}{V^{(i)}(x_i,\sigma)}\frac{t}{\psi_j^{(i)}(\sigma)}\left(\frac{t}{\phi^{(i)}(\sigma)}+1\right)^{\alpha_{2,0}}\label{hk-up-tr}\\
\le&\frac{2^{\alpha_{2,0}}c'_6C_{VD,0}^2C_{DUE,0}}{V^{(i)}(x_i,\phi^{-i}(t))}\left(1+\frac{d^{(i)}(x_i,y_i)}{\phi^{-i}(t)}\right)^{\alpha_{2,0}}\exp\left\{\left(\left(c'_5+\frac{c'_7}{2}\right)\vee 2\right)\frac{t}{\phi^{(i)}(\sigma)}-\frac{c'_8}{2}\frac{d^{(i)}(x_i,y_i)}{\sigma}\right\}+\notag\\
&\quad\quad\quad+\frac{2C_{J,0}c'_9}{V^{(i)}(x_i,\sigma)}\left(\frac{t}{\phi^{(i)}(\sigma)}+1\right)^{\alpha_{2,0}+1}.\notag
\end{align}
The lemma follows by taking $c_6=2^{\alpha_{2,0}}c'_6C_{VD,0}^2C_{DUE,0},$ $c_7=c'_5+c'_7/2,$ $c_8=c'_8/2$ and $c_9=2C_{J,0}c'_9$.
\end{proof}

Note that $c_7\ge 2$ and $c_8\le\frac{1}{2}$ according to the explicit formulae.

\begin{lemma}\cite[Lemma 3.4 covered]{bkkl}\label{3.4}
Assume that $\pi:\mathbb{R}_+\times\mathbb{R}_+\times I\to\mathbb{R}_+$ and $\{T_i\}_{i\in I}$ (where each $T_i\in(0,\infty]$) satisfy the following conditions:

(i) $\pi(\cdot,t;i)$ is non-decreasing, and $\pi(r,\cdot;i)$ is non-increasing;

(ii) For each $b>0$, there exists $C^*_\pi(b)<\infty$ such that for all $i\in I$ and $0<t<T_i$, $\pi(b\phi^{-i}(t),t;i)\le C^*_\pi(b)$;

(iii) There exist $\eta\in(0,\beta_{2,0}]$ and $a_1,a_2>0$ such that for all $i\in I$, $x_i\in M^{(i)},r>0$ and $0<t<T_i$,
\begin{equation}\label{eta-tail}
P^{(i)}_t1_{B_i(x_i,r)^c}(x)\le a_2\left(\frac{\psi_j^{-i}(t)}{r}\right)^\eta+a_2\exp\{-a_1\pi(r,t;i)\}.
\end{equation}
Then under the same assumptions as in Proposition \ref{found+}, there exist a non-negative integer $k_\eta$ and a constant $c_{10}=c_{10}(\eta,a_1,a_2;\pi)>0$ such that for all $i\in I$, $x_i,y_i\in M^{(i)}$ and $0<t<T_i$,
$$p^{(i)}_t(x_i,y_i)\le\frac{c_{10}t}{V^{(i)}(x_i,y_i)\psi_j^{(i)}(x_i,y_i)}+\frac{c_{10}}{V^{(i)}(x_i,\phi^{-i}(t))}\left(1+\frac{d^{(i)}(x_i,y_i)}{\phi^{-i}(t)}\right)^{\alpha_{2,0}}\exp\left\{-a_1k_\eta\pi\left(\frac{d^{(i)}(x_i,y_i)}{16k_\eta},t;i\right)\right\}.$$
\end{lemma}

\begin{proof}
By \cite[Formulae (3.7) and (3.8)]{bkkl}, for every $B_i=B_i(x_i,r)$, for $\mu^{(i)}$-a.e. $z_i\in\frac{1}{4}B_i$ and all $0<t<T_i$,
$$1-P^{(i,B_i)}_t1_{B_i}(z_i)\le 2a_2\left(\frac{\psi_j^{-i}(2t)}{r/4}\right)^\eta+2a_2\exp\{-a_1\pi(r/4,2t)\}+\frac{2c'_5t}{\psi_j^{(i)}(r)}.$$
Let $k_\eta=\lceil\eta^{-1}(\beta_{2,0}+2\alpha_{2,0})\rceil$. Then by \cite[Formulae (3.9) to (3.12)]{bkkl}, for $\sigma_i=d^{(i)}(x_i,y_i)/(4k_\eta)>\phi^{-i}(t)$,
$$p^{(i,\sigma_i)}_t(x_i,y_i)\le\frac{c'_{10}(a_1,a_2;\pi)}{V^{(i)}(x_i,\phi^{-i}(t))}\left(1+\frac{d^{(i)}(x_i,y_i)}{\phi^{-i}(t)}\right)^{\alpha_{2,0}}\left(\left(\frac{\psi_j^{-i}(t)}{d^{(i)}(x_i,y_i)}\right)^{\beta_{2,0}+2\alpha_{2,0}}+\exp\left\{-a_1k_\eta\pi\left(\frac{d^{(i)}(x_i,y_i)}{16k_\eta},t\right)\right\}\right),$$
where
$$c'_{10}(\eta,a_1,a_2;\pi)=3^{k_\eta}C_{VD,0}C_{DUE,0}\exp(2c'_5a_1k_\eta C^*_\pi(1/4))\left\{(4k_\eta)^{\beta_{2,0}+2\alpha_{2,0}}\left((2^{2\beta_{1,0}+1}a_2)^{k_\eta}+(c'_5/C_{\psi,0})^{k_\eta}\right)+(2a_2)^{k_\eta}\right\}.$$
Thus parallel to (\ref{hk-up-tr}), we obtain
\begin{align*}
p^{(i)}_t(x_i,y_i)\le&\frac{c'_{10}(\eta,a_1,a_2;\pi)}{V^{(i)}(x_i,\phi^{-i}(t))}\left(1+\frac{d^{(i)}(x_i,y_i)}{\phi^{-i}(t)}\right)^{\alpha_{2,0}}\exp\left\{-a_1k_\eta\pi\left(\frac{d^{(i)}(x_i,y_i)}{16k_\eta},t\right)\right\}+\\
&+\frac{c'_{10}(\eta,a_1,a_2;\pi)}{V^{(i)}(x_i,\phi^{-i}(t))}\left(\frac{2d^{(i)}(x_i,y_i)}{\phi^{-i}(t)}\right)^{\alpha_{2,0}}\left(\frac{\psi_j^{-i}(t)}{d^{(i)}(x_i,y_i)}\right)^{\beta_{2,0}+2\alpha_{2,0}}+\frac{2C_{J,0}C_9t}{V^{(i)}(x_i,\sigma_i)\psi_j^{(i)}(\sigma_i)}\left(\frac{t}{\phi^{(i)}(\sigma_i)}+1\right)^{\alpha_{2,0}}\\
\le&\frac{c'_{10}(\eta,a_1,a_2;\pi)}{V^{(i)}(x_i,\phi^{-i}(t))}\left(1+\frac{d^{(i)}(x_i,y_i)}{\phi^{-i}(t)}\right)^{\alpha_{2,0}}\exp\left\{-a_1k_\eta\pi\left(\frac{d^{(i)}(x_i,y_i)}{16k_\eta},t\right)\right\}+\\
&+\frac{c'_{10}(\eta,a_1,a_2;\pi)C_{VD,0}C_{\psi,0}(2C_{\psi,0})^{\alpha_{2,0}+1}t}{V^{(i)}(x_i,y_i)\psi_j^{(i)}(x_i,y_i)}+\frac{2^{\alpha_{2,0}+1}C_{VD,0}C_{\psi,0}C_{J,0}C_9(4k_\eta)^{\alpha_{2,0}+\beta_{2,0}}t}{V^{(i)}(x_i,y_i)\psi_j^{(i)}(x_i,y_i)}.
\end{align*}
Hence the proof is completed by letting
$$c_{10}(\eta,a_1,a_2;\pi):=2^{\alpha_{2,0}+1}C_{VD,0}C_{\psi,0}^{\alpha_{2,0}+2}c'_{10}(\eta,a_1,a_2;\pi)+2^{\alpha_{2,0}+1}C_{VD,0}C_{\psi,0}C_{J,0}C_9(4k_\eta)^{\alpha_{2,0}+\beta_{2,0}},$$
which is independent of $i\in I$.
\end{proof}

For any $i\in I$, let $r'_{*,i}:=r_{*,i}\wedge\overline{R}_i$. If $r'_{*,i}\in(0,\infty)$, define $\lambda_{T,i}:=(\phi^{-i}(T)/r'_{*,i})\vee 1$ for all $T\in(0,\infty)$; while if $r'_{*,i}=\infty$, define $\lambda_{T,i}:=1$ for all $T\in(0,\infty]$.

\begin{lemma}\cite[Lemmas 3.8 and 3.9 covered]{bkkl}\label{3.8}
Under the same assumptions as in Proposition \ref{found+}, for all $i\in I$, $T\in(0,\infty)$ (if $r'_{*,i}<\infty$) or $T\in(0,\infty]$ (if $r'_{*,i}=\infty$), $r>0$ and $0<t<T$,
\begin{equation}\label{minimax}
\phi^{(i)}_*(r,t):=\sup\limits_{\sigma>0}\left(\frac{r}{\sigma}-\frac{t}{\phi^{(i)}(\sigma)}\right)=\sup\limits_{\sigma_1(r,t;T,i)\le\sigma\le 2\phi^{-i}(t)}\left(\frac{r}{\sigma}-\frac{t}{\phi^{(i)}(\sigma)}\right)\ge\frac{r}{2\phi^{-i}(t)},
\end{equation}
where
$$\sigma_1(r,t;T,i):=\left(\lambda_{T,i}^{-1}\phi^{-i}(t)\right)^{\frac{\beta_{*,0}}{\beta_{*,0}-1}}(C_{\phi,0}r)^{-\frac{1}{\beta_{*,0}-1}}.$$
Further, given any $\delta_1,\delta_2,\delta_3>0$, for all $0<t\le T$, $\phi^{(i)}_*(\delta_1r,\delta_2t)\le c_{11}(\delta_1,\delta_2;T,i)\phi^{(i)}_*(r,t)$ if $r\ge 2\phi^{-i}(t)$, while $\phi^{(i)}_*(r,t)\le c_{12}(\delta_3;T,i)$ if $r\le\delta_3\phi^{-i}(t)$, where
\begin{align*}
c_{11}(\delta_1,\delta_2;T,i):=\left(\frac{C_{\phi,0}}{\delta_2}\right)^{\frac{1}{\beta_{*,0}-1}}(\lambda_{T,i}(\delta_1\vee\delta_2))^{\frac{\beta_{*,0}}{\beta_{*,0}-1}},\quad c_{12}(\delta_3;T,i):=2^{\frac{\beta_{2,0}-\beta_{*,0}}{\beta_{*,0}-1}}(\beta_{*,0}-1)C_{\phi,0}^{\frac{2}{\beta_{*,0}-1}}(\lambda_{T,i}(\delta_3\vee 2))^{\frac{\beta_{*,0}}{\beta_{*,0}-1}}.
\end{align*}
\end{lemma}

\begin{proof}
These inequalities completely follow from the proofs in \cite{bkkl}.
\end{proof}

For simplicity, we denote $\beta_{*,1}:=1/(\beta_{*,0}-1)$. Note that if $r'_{*,i}<\infty$, for each $T\le T_i:=\phi^{(i)}(r'_{*,i})$,
$$\lambda_{T,i}=\frac{\phi^{-1}(T)}{r'_{*,i}}\vee 1\le\frac{\phi^{-i}(T_i)}{r'_{*,i}}\vee 1=1.$$
The same is already given by the definition if $r'_{*,i}=\infty$. Therefore, $c_{11}(\delta_1,\delta_2;T,i)=(C_{\phi,0}/\delta_2)^{\beta_{*,1}}(\delta_1\vee\delta_2)^{\beta_{*,1}+1}=:c_{11}(\delta_1,\delta_2)$ and $c_{12}(\delta_3;T,i)=2^{\beta_{*,1}(\beta_{2,0}-\beta_{*,0})}\beta_{*,1}^{-1}C_{\phi,0}^{2\beta_{*,1}}(\delta_3\vee 2)^{\beta_{*,0}\beta_{*,1}}=:c_{12}(\delta_3)$ for $T\le T_i$, independent of $i$.

\begin{lemma}\cite[Inequality (3.25) covered]{bkkl}
Under the same assumptions as in Proposition \ref{found+}, there exist $c_{13},c_{14}>0$ such that for all $i\in I$, $x_i\in M^{(i)}$, $r>0$ and $0<t<T_i$,
\begin{equation}\label{3.25'}
P_t1_{B_i(x_i,r)^c}(x_i)\le c_{13}\left(\frac{\psi_j^{-i}(t)}{r}\right)^{\beta_{1,0}/2}+c_{13}\exp(-c_{14}\phi_*^{(i)}(r,t)).
\end{equation}
\end{lemma}

\begin{proof}
We split the proof into three cases:

(1) For all $0<r\le\frac{4}{c_8}\phi^{-i}(c_7t)$, by Lemma \ref{3.8},
$$\phi_*^{(i)}(r,t)\le\phi_*^{(i)}\left(\frac{4}{c_8}\phi^{-i}(c_7t),t\right)\le c_{11}(1,c_7^{-1};T_i,i)c_{12}(4/c_8;T_i/c_7,i)=c_{11}(1,c_7^{-1})c_{12}(4/c_8).$$
Therefore, for each $c>0$,
$$P_t^{(i)}1_{B_i(x_i,t)^c}\le 1\le\exp(cc_{11}(1,c_7^{-1})c_{12}(4/c_8))\cdot \exp(-c\phi_*^{(i)}(r,t))=:\chi(c)\exp(-c\phi_*^{(i)}(r,t)).$$

(2) For all $\frac{4}{c_8}\phi^{-i}(c_7t)<r\le\frac{4}{c_8}\frac{\phi^{-i}(c_7t)^{1+\theta}}{\psi_j^{-i}(c_7t)^\theta}$, where
$$\theta=\frac{\beta_{1,0}(\beta_{*,0}-1)}{\beta_{*,0}(2\alpha_{2,0}+\beta_{1,0})+(\beta_{1,0}+2\beta_{2,1}+2\alpha_{2,0}\beta_{2,1})}\in(0,1):$$
For each $\sigma\le 2\phi^{-i}(c_7t)$ and $\alpha>0$, let $\sigma_n=2^{n\alpha}\sigma$. Then (cf. \cite[Inequality (3.26)]{bkkl}),
\begin{align}
\frac{c_7t}{\phi^{(i)}(\sigma_n)}-\frac{c_82^nr}{\sigma_n}\le&-\left((C_{\phi,0}c_7)^{(\beta_{*,1}+1)/\beta_{1,1}}c_{11}(c_8/2,c_7)\right)^{-1}\phi_*^{(i)}(r,t)-\left(\frac{c_8}{8}(C_{\phi,0}c_7)^{-1/\beta_{1,1}}\right)2^{n(1-\alpha)}\frac{r}{\phi^{-i}(t)}\notag\\
=:&-c'_{11}\phi_*^{(i)}(r,t)-c'_{12}\frac{2^{n(1-\alpha)}r}{\phi^{-i}(t)}.\label{3.26'}
\end{align}
For every $\frac{\alpha_{2,0}}{\alpha_{2,0}+\beta_{1,0}}<\alpha<1$, by Lemma \ref{3.2} we obtain for all $t\le T_i$ and $y_i\in B_i(x_i,2^{n+1}r)\setminus B_i(x_i,2^nr)$ (cf. \cite[Inequality (3.27)]{bkkl}) that
$$p_t^{(i)}(x_i,y_i)\le\frac{c_6C_*(\alpha_{2,0},c'_{12}/2,1-\alpha)}{V(x,\phi^{-i}(t))}\exp\left(-c'_{11}\phi_*^{(i)}(r,t)-\frac{c'_{12}}{2}\frac{2^{n(1-\alpha)}r}{\phi^{-i}(t)}\right)+\frac{c_9t}{V^{(i)}(x_i,\sigma_n)\psi_j^{(i)}(\sigma_n)}\left(1+\frac{t}{\sigma_n}\right)^{\alpha_{2,0}+1},$$
and therefore (cf. \cite[the paragraphs between (3.27) and (3.29)]{bkkl})
$$P_t^{(i)}1_{B_i(x_i,r)^c}(x_i)\le c'_{13}\exp(-c'_{11}\phi_*^{(i)}(r,t))+c'_{14}\left(\psi_j^{-i}(c_7t)/r\right)^{\beta_{1,0}/2},$$
where
\begin{align*}
c'_{13}=&\frac{c_6C_{VD,0}}{1-2^{-\alpha_{2,0}}}C_*\left(\alpha_{2,0},\frac{c'_{12}}{2},1-\alpha\right)C_*\left(\alpha_{2,0},\frac{c'_{12}}{4},1\right)\exp\left(\alpha_{2,0}c_*\left(1-\alpha,\frac{c'_{12}}{4\alpha_{2,0}}\right)\right),\\
c'_{14}=&2^{\beta_{2,0}\beta_{*,1}}C_{VD,0}C_{\psi,0}^2C_{\phi,0}\frac{c_9(1+c_7^{-1})^{\alpha_{2,0}+2}}{1-2^{\alpha_{2,0}-\alpha(\alpha_{2,0}+\beta_{1,0})}}\left(\frac{c_8}{4}\right)^{\frac{\beta_{1,0}}{2\theta}}\left(c_7C_{\phi,0}^{1+\frac{\beta_{*,0}}{\beta_{1,1}}}\right)^{\beta_{*,1}(\alpha_{2,0}+\beta_{1,0}+\beta_{2,0}+\beta_{2,1}(1+\alpha_{2,0}))}.
\end{align*}
In particular,
$$P_t^{(i)}1_{B_i(x_i,r)^c}(x_i)\le c'_{13}\exp(-c'_{11}\phi_*^{(i)}(r,t))+c'_{14}\sqrt{c_7C_{\psi,0}}\left(\psi_j^{-i}(t)/r\right)^{\beta_{1,0}/2}.$$

(3) For all $r>\frac{4}{c_8}\frac{\phi^{-i}(c_7t)^{1+\theta}}{\psi_j^{-i}(c_7t)^\theta}$, we simplify \cite[the proof of Lemma 3.6]{bkkl}.

Let $N=\lceil 2\frac{\alpha_{2,0}}{\beta_{1,0}}\rceil+2$ so that $\eta=\beta_{1,0}-\frac{\beta_{1,0}+\alpha_{2,0}}{N}\ge\frac{\beta_{1,0}}{2}$. Following the proof of \cite[Inequality (3.14)]{bkkl}, by taking $\alpha=(2\alpha_{2,0}+\beta_{1,0})/(2\alpha_{2,0}+2\beta_{1,0})$, we have for all $i\in I$, $x_i\in M^{(i)}$ and $r',t>0$,
\begin{equation}\label{3.14'}
P^{(i)}_t1_{B_i(x_i,r')^c}(x)\le c'_{15}\left(\frac{\psi_j^{-i}(t)}{r'}\right)^\eta+c'_{16}\exp\left(-\frac{c_8}{2}\left(\frac{r'}{\phi^{-i}(t)}\right)^{1/N}\right),
\end{equation}
where
\begin{align*}
c'_{15}=&e^{\frac{c_8}{2}}\vee\frac{2^{2\alpha_{2,0}+1}c_9C_{VD,0}C_{\psi,0}}{1-2^{-\beta_{1,0}/2}},\\
c'_{16}=&\frac{2^{3\alpha_{2,0}+2}c_6C_{VD,0}}{1-2^{-(2\alpha_{2,0}+1)}}\exp\left\{c_7+(2\alpha_{2,0}+1)c_*\left(\frac{\beta_{1,0}}{2(\alpha_{2,0}+\beta_{1,0})},\frac{c_8}{2(2\alpha_{2,0}+1)}\right)\right\}C_*\left(2\alpha_{2,0}+1,\frac{c_8}{4},\frac{1}{N}\right).
\end{align*}
Since $r>\frac{4}{c_8}\frac{\phi^{-i}(c_7t)^{1+\theta}}{\psi_j^{-i}(c_7t)^\theta}$, there exists $\theta_0=\theta_0(r,t)>\theta$ such that $r=\frac{4}{c_8}\frac{\phi^{-i}(c_7t)^{1+\theta_0}}{\psi_j^{-i}(c_7t)^{\theta_0}}$. Therefore,
\begin{align*}
\exp&\left\{-\frac{c_8}{2}\left(\frac{r}{\phi^{-i}(t)}\right)^{\frac{1}{N}}\right\}\le\exp\left\{-\frac{c_8}{2(c_7C_{\phi,0})^{1/\beta_{1,1}}}\left(\frac{r}{\phi^{-i}(c_7t)}\right)^{\frac{1}{N}}\right\}=\exp\left\{-\frac{c_8}{2(c_7C_{\phi,0})^{1/\beta_{1,1}}}\left(\frac{\phi^{-i}(c_7t)}{\psi_j^{-i}(c_7t)}\right)^{\frac{\theta_0}{N}}\right\}\\
\le&C_*\left(\eta+\frac{\eta}{\theta_0},\frac{c_8}{2(c_7C_{\phi,0})^{1/\beta_{1,1}}},\frac{1}{N}\right)\left(\frac{\phi^{-i}(c_7t)}{\psi_j^{-i}(c_7t)}\right)^{-\eta(1+\theta_0)}\le C_*\left(\eta+\frac{\eta}{\theta},\frac{c_8}{2(c_7C_{\phi,0})^{1/\beta_{1,1}}},\frac{1}{N}\right)\left(\frac{4}{c_8}\frac{\psi_j^{-i}(c_7t)}{r}\right)^\eta\\
\le&(4c_8^{-1})^\eta(c_7C_{\psi,0})^{\eta/\beta_{1,0}}C_*\left(\eta+\frac{\eta}{\theta},\frac{c_8}{2(c_7C_{\phi,0})^{1/\beta_{1,1}}},\frac{1}{N}\right)\left(\frac{\psi_j^{-i}(t)}{r}\right)^\eta=:c'_{17}\left(\frac{\psi_j^{-i}(t)}{r}\right)^\eta.
\end{align*}
Thus by (\ref{3.14'}) with $r'=r$, we obtain
$$P^{(i)}_t1_{B_i(x_i,r)^c}(x)\le(c'_{15}+c'_{17})\left(\frac{\psi_j^{-i}(t)}{r}\right)^\eta\le(c'_{15}+c'_{16}c'_{17})\left(\frac{\psi_j^{-i}(t)}{r}\right)^{\beta_{1,0}/2}.$$

Combining the cases (1), (2) and (3), we see (\ref{3.25'}) holds with $c_{13}=\chi(c'_{11})\vee c'_{14}\sqrt{c_7C_{\psi,0}}\vee c'_{13}\vee(c'_{15}+c'_{16}c'_{17})$ and $c_{14}=c'_{11}$.
\end{proof}

\begin{proof}[Proof of Proposition \ref{found+}]
Parallel to the proof of \cite[Theorem 2.13 (i)]{bkkl}. Let $k=k_{\beta_{1,0}/2}$. For all $i\in I$, $t<T_i$, if $d^{(i)}(x_i,y_i)\le 32k\phi^{-i}(t)$, then $\phi_*^{(i)}(d^{(i)}(x_i,y_i),t)\le c_{12}(32k;T_i,i)=c_{12}(32k)$. By the uniform DUE condition, we obtain (for all $c>0$)
$$p_t^{(i)}(x_i,y_i)\le\frac{C_{DUE,0}}{V^{(i)}(x_i,\phi^{-i}(t))}\le\frac{C_{DUE,0}\exp\{cc_{12}(32k)\}}{V^{(i)}(x_i,\phi^{-i}(t))}\exp(-c\phi_*^{(i)}(d^{(i)}(x_i,y_i),t)).$$

If $d^{(i)}(x_i,y_i)>32k\phi^{-i}(t)$, let $\pi(r,t;i)=\phi_*^{(i)}(r,t)$. By inserting (\ref{3.25'}) and (\ref{3.26'}) (with $n=0$) into Lemma \ref{3.4}, we obtain (here $c_{10}:=c_{10}(\beta_{1,0}/2,c_{14},c_{13};\pi)$)
\begin{align*}
p^{(i)}_t(x_i,y_i)-&\frac{c_{10}t}{V^{(i)}(x_i,y_i)\psi_j^{(i)}(x_i,y_i)}\le\frac{c_{10}}{V^{(i)}(x_i,\phi^{-i}(t))}\left(1+\frac{d^{(i)}(x_i,y_i)}{\phi^{-i}(t)}\right)^{\alpha_{2,0}}\exp\left\{-c_{14}k\phi_*^{(i)}\left(\frac{d^{(i)}(x,y)}{16k},t\right)\right\}\\
\le&\frac{c_{10}}{V^{(i)}(x_i,\phi^{-i}(t))}\left(1+\frac{d^{(i)}(x_i,y_i)}{\phi^{-i}(t)}\right)^{\alpha_{2,0}}\exp\left\{-\frac{c_{14}k}{c_{11}(16kc_8,c_7;T,i)}\phi_*^{(i)}(c_8d^{(i)}(x,y),c_7t)\right\}\\
\le&\frac{c_{10}}{V^{(i)}(x_i,\phi^{-i}(t))}\left(1+\frac{d^{(i)}(x_i,y_i)}{\phi^{-i}(t)}\right)^{\alpha_{2,0}}\exp\left\{-\frac{c_{14}k}{c_{11}(16kc_8,c_7)}\left(c'_{11}\phi_*^{(i)}(d^{(i)}(x_i,y_i),t)+c'_{12}\frac{d^{(i)}(x_i,y_i)}{\phi^{-i}(t)}\right)\right\}\\
\le&\frac{c_{10}C_*(\alpha_{2,0},\frac{c_{14}c'_{12}}{c_{11}(16kc_8,c_7)},1)}{V^{(i)}(x_i,\phi^{-i}(t))}\exp\left\{-\frac{c_{14}c'_{11}k}{c_{11}(16kc_8,c_7)}\phi_*^{(i)}(d^{(i)}(x_i,y_i),t)\right\}
\end{align*}
Combining the two cases, we see (\ref{uhk}) holds with constants
$$c=\frac{c_{14}c'_{11}k_{\beta_{1,0}/2}}{c_{11}(16c_8k_{\beta_{1,0}/2},c_7)},\quad C=\left(C_{DUE,0}\exp\left\{cc_{12}\left(32k_{\beta_{1,0}/2}\right)\right\}\right)\vee c_{10}\vee\left(c_{10}C_*\left(\alpha_{2,0},\frac{c_{14}c'_{12}}{c_{11}(16c_8k_{\beta_{1,0}/2},c_7)},1\right)\right).$$
Hence we have proved the desired result.
\end{proof}

\subsection{The two-sided heat kernel for jump processes}\label{else}

Now we prove the best heat kernel estimate under the $\rho$-chain condition. For simplicity, we denote
$$V_{\varepsilon}(x,r):=\mu(B_{d_\varepsilon}(x,r))\quad\mbox{and}\quad\overline{R}_\varepsilon:=diam(M,d_\varepsilon).$$

\begin{lemma}\label{key}
Assume that the $\rho$-chain condition holds, and that there exists a regular strongly local Dirichlet form $(\mathcal{E}_c,\mathcal{F}_c)$ satisfying $\mathrm{PHI}(\psi_c)$. Then, $\mathrm{VD}$ under the metric $d$ implies $\mathrm{VD}$ under the metric $d_\varepsilon$, uniformly on $\varepsilon>0$. Assume that $(\mathcal{E}_j,\mathcal{F}_j)$ is a regular jump-type Dirichlet form whose jump kernel is $J$. Then, uniformly on $\varepsilon>0$, conditions $\mathrm{DUE}(\Phi)$, $\mathrm{E}(\Phi)$ and $\mathrm{J}(\psi_j)$ (under the metric $d$) imply respectively $\mathrm{DUE}(\Phi\circ F_{\rho,\varepsilon})$, $\mathrm{E}(\Phi\circ F_{\rho,\varepsilon})$ and $\mathrm{J}(\psi_j\circ F_{\rho,\varepsilon})$ (under the metric $d_\varepsilon$). 
\end{lemma}

\begin{proof}
It is clear by Lemma \ref{beta-uni} and Corollary \ref{beta-uni'} that the scale functions $\Phi\circ F_{\rho,\varepsilon}$ or $\psi_j\circ F_{\rho,\varepsilon}$ are uniform on $\varepsilon$.

By Lemma \ref{rho-lem} and Corollary \ref{rho-cor}, there exists a constant $C\ge 1$ such that for any $x\in M$ and $\varepsilon,r>0$,
\begin{equation}\label{reg-rho}
B_d(x,C^{-1}F_{\rho,\varepsilon}(r))\subset B_{d_\varepsilon}(x,r)\subset B_d(x,CF_{\rho,\varepsilon}(r)),\mbox{ and }B_{d_\varepsilon}(x,C^{-1}F_{\rho,\varepsilon}^{-1}(r))\subset B_d(x,r)\subset B_{d_\varepsilon}(x,CF_{\rho,\varepsilon}^{-1}(r)).
\end{equation}
By the first relation we have $2\overline{R}>C^{-1}F_{\rho,\varepsilon}(\frac{1}{2}\overline{R}_\varepsilon)$, while by the third we have $2\overline{R}_\varepsilon>C^{-1}F_{\rho,\varepsilon}^{-1}(\frac{1}{2}\overline{R})$, if finite.

Therefore, uniformly on $\varepsilon>0$, by (VD), for all $x\in M$ and $R>r>0$,
$$\frac{V_\varepsilon(x,R)}{V_\varepsilon(x,r)}=\frac{\mu(B_{d_\varepsilon}(x,R))}{\mu(B_{d_\varepsilon}(x,r))}\le\frac{\mu(B_d(x,CF_{\rho,\varepsilon}(R)))}{\mu(B_d(x,C^{-1}F_{\rho,\varepsilon}(r)))}\le C\left(\frac{CF_{\rho,\varepsilon}(R)}{C^{-1}F_{\rho,\varepsilon}(r)}\right)^{\alpha_2}\le C'\left(\frac{R}{r}\right)^{\alpha_2},$$
so that the spaces $(M,d_\varepsilon,\mu)$ satisfy (VD) uniformly.%; while on the other hand by (RVD) for all $0<r<R<C^{-1}F_{\rho,\varepsilon}(R)$,
%$$\frac{V_\varepsilon(x,R)}{V_\varepsilon(x,r)}\ge\frac{\mu(B_d(x,C^{-1}F_{\rho,\varepsilon}(R)))}{\mu(B_d(x,CF_{\rho,\varepsilon}(r)))}\ge C^{-1}\left(\frac{C^{-1}F_{\rho,\varepsilon}(R)}{CF_{\rho,\varepsilon}(r)}\right)^{\alpha_1}\ge (C')^{-1}\left(\frac{R}{r}\right)^{\alpha_1/\gamma_2},$$
%so that the spaces $(M,d_\varepsilon,\mu)$ satisfy (RVD) uniformly.

Further, by (\ref{reg-rho}), the topology induced by $d_\varepsilon$ is identical to the topology induced by $d$. Therefore, the Dirichlet space $(M,d_\varepsilon,\mu,\mathcal{E},\mathcal{F})$ is regular if and only if $(M,d,\mu,\mathcal{E},\mathcal{F})$ is regular.

Now we turn to the conditions related to the jump kernel. Uniformly on $\varepsilon>0$,

(1) for every Dirichlet form satisfying $\mathrm{DUE}(\Phi)$, for all $x,y\in M$ and $t<\Phi(\overline{R})$,
\begin{equation}\label{due-0}
p_t(x,y)\le\frac{C}{V(x,\Phi^{-1}(t))}\le\frac{C}{V_\varepsilon(x,C^{-1}F_{\rho,\varepsilon}^{-1}(\Phi^{-1}(t)))}\le\frac{C}{V_\varepsilon(x,C^{-2}F_{\rho,\varepsilon}^{-1}\Phi^{-1}(t))}\le\frac{C'}{V_\varepsilon(x,(\Phi\circ F_{\rho,\varepsilon})^{-1}(t))}.
\end{equation}
If $\overline{R}=\infty$, then this is exactly the condition $\mathrm{DUE}(\Phi\circ F_{\rho,\varepsilon})$ under the metric $d_\varepsilon$. If $\overline{R}<\infty$, note that
$$k=\left\lceil\frac{\Phi\circ F_{\rho,\varepsilon}(\overline{R}_\varepsilon)}{\Phi(\overline{R})}\right\rceil\le\left\lceil\frac{2^{\beta_2(\Phi)}\Phi(2C\overline{R})}{\Phi(\overline{R})}\right\rceil\le\left\lceil C(4C)^{\beta_2(\Phi)}\right\rceil=:k_0$$
uniformly on $\varepsilon$. Same as in Remark \ref{release-T}, by (\ref{due-0}) we have for all $t<\Phi\circ F_{\rho,\varepsilon}(\overline{R}_\varepsilon)$ that
\begin{align*}
p_t(x,y)=&P_{(k_0-1)t/k_0}p_{t/k_0}(x,\cdot)(y)\le\|p_{t/k_0}(x,\cdot)\|_{L^\infty}\le\frac{C'}{V_\varepsilon(x,(\Phi\circ F_{\rho,\varepsilon})^{-1}(t/k_0))}\le\frac{C''(k_0)}{V_\varepsilon(x,(\Phi\circ F_{\rho,\varepsilon})^{-1}(t))}.
\end{align*}
Therefore, condition $\mathrm{DUE}(\Phi\circ F_{\rho,\varepsilon})$ holds under the metric $d_\varepsilon$;

(2) for every Dirichlet form satisfying $\mathrm{E}(\Phi)$, if $CF_{\rho,\varepsilon}(r)<\kappa_E\overline{R}$, then on one hand,
$$G^{B_{d_\varepsilon}(x,r)}1_{B_{d_\varepsilon}(x,r)}\le G^{B_d(x,CF_{\rho,\varepsilon}(r))}1_{B_d(x,CF_{\rho,\varepsilon}(r))}\le C\Phi(CF_{\rho,\varepsilon}(r))\le C'\Phi\circ F_{\rho,\varepsilon}(r);$$
on the other hand,
\begin{align*}
G^{B_{d_\varepsilon}(x,r)}1_{B_{d_\varepsilon}(x,r)}\ge&G^{B_d(x,C^{-1}F_{\rho,\varepsilon}(r))}1_{B_d(x,C^{-1}F_{\rho,\varepsilon}(r))}\ge C^{-1}\Phi(C^{-1}F_{\rho,\varepsilon}(r))1_{B_d(x,(4C)^{-1}F_{\rho,\varepsilon}(r))}\\
\ge&(C')^{-1}\Phi\circ F_{\rho,\varepsilon}(r)1_{B_{d_\varepsilon}(x,(C'')^{-1}r)}.
\end{align*}
By the standard covering argument, for each $\xi\in B_{d_\varepsilon}(x,r/4)$ we have
$$G^{B_{d_\varepsilon}(x,r)}1_{B_{d_\varepsilon}(x,r)}(\xi)\ge G^{B_{d_\varepsilon}(\xi,\frac{3}{4}r)}1_{B_{d_\varepsilon}(\xi,\frac{3}{4}r)}(\xi)\ge(C')^{-1}\Phi\circ F_{\rho,\varepsilon}\left(\frac{3}{4}r\right)1_{B_{d_\varepsilon}(\xi,(C'')^{-1}\frac{3}{4}r)}(\xi)\ge(\tilde{C})^{-1}\Phi\circ F_{\rho,\varepsilon}(r).$$
Thus it suffices to find $0<\kappa<1$ such that $\kappa\overline{R}_\varepsilon\le F_{\rho,\varepsilon}^{-1}(C^{-1}\kappa_E\overline{R})$ for every $\varepsilon>0$, if $\overline{R}<\infty$. Since
$$\frac{\overline{R}_\varepsilon}{F_{\rho,\varepsilon}^{-1}(C^{-1}\kappa_E\overline{R})}\le\frac{2F_{\rho,\varepsilon}^{-1}(2C\overline{R})}{F_{\rho,\varepsilon}^{-1}(C^{-1}\kappa_E\overline{R})}\le 4C^3\kappa_E^{-1},$$
we can take $\kappa:=\frac{1}{4}C^{-3}\kappa_E$, and condition $\mathrm{E}(\Phi\circ F_{\rho,\varepsilon})$ holds under the metric $d_\varepsilon$;

(3) for every jump kernel $J$ satisfying $\mathrm{J}(\psi_j)$, on one hand,
\begin{align*}
J(x,y)\le&\frac{C}{V(x,d(x,y))\psi_j(d(x,y))}\le\frac{C}{V_\varepsilon(x,C^{-1}F_{\rho,\varepsilon}^{-1}(d(x,y)))\psi_j(C^{-1}F_{\rho,\varepsilon}(d_\varepsilon(x,y)))}\\
\le&\frac{C'}{V_\varepsilon(x,C^{-2}d_\varepsilon(x,y))\psi_j(F_{\rho,\varepsilon}(d_\varepsilon(x,y)))}\le\frac{C''}{V_\varepsilon(x,d_\varepsilon(x,y))\psi_j\circ F_{\rho,\varepsilon}(d_\varepsilon(x,y))};
\end{align*}
on the other hand, similarly
\begin{align*}
J(x,y)\ge&\frac{C^{-1}}{V(x,d(x,y))\psi_j(d(x,y))}\ge\frac{C^{-1}}{V_\varepsilon(x,CF_{\rho,\varepsilon}^{-1}(d(x,y)))\psi_j(CF_{\rho,\varepsilon}(d_\varepsilon(x,y)))}\\
\ge&\frac{(C')^{-1}}{V_\varepsilon(x,C^2d_\varepsilon(x,y))\psi_j(F_{\rho,\varepsilon}(d_\varepsilon(x,y)))}\ge\frac{(C'')^{-1}}{V_\varepsilon(x,d_\varepsilon(x,y))\psi_j\circ F_{\rho,\varepsilon}(d_\varepsilon(x,y))}.
\end{align*}
Combining these two inequalities, we see $\mathrm{J}(\psi_j\circ F_{\rho,\varepsilon})$ holds uniformly on $\varepsilon>0$.
\end{proof}

%$$\tilde{\Phi}(r)=\frac{\psi_c\circ F_{\rho,\varepsilon}(r)}{\int_0^r\frac{d\psi_c\circ F_{\rho,\varepsilon}(s)}{\psi_j\circ F_{\rho,\varepsilon}(s)}}\xlongequal{\sigma=F_{\rho,\varepsilon}(s)}\frac{\psi_c(F_{\rho,\varepsilon}(r))}{\int_0^r\frac{d\psi_c(\sigma)}{\psi_j(\sigma)}}=\Phi(F_{\rho,\varepsilon}(r)).$$

\begin{proof}[Proof of Theorem \ref{thm2}]
(1) We first prove the lower bound. Let us begin by explaining how the subordination techniques work in the bounded case. We first replace the L\'evy measure $\nu(dt)=\frac{dt}{t\psi_j\circ\psi_c^{-1}(t)}$ by $\nu_{\overline{R}}(dt)=\frac{1_{t<\overline{R}}dt}{t\psi_j\circ\psi_c^{-1}(t)}$. By \cite[Lemmas 3.1 and 3.3]{lm}, the subordinated Dirichlet form $(\mathcal{E}_{j,0},\mathcal{F}_{j,0})$ is still jump-type with a jump kernel satisfying $\mathrm{J}(\psi_j)$. Its $\mathrm{LLE}(\Phi)$ can be proved in the same way as in \cite[Lemma 4.5]{bkkl}, which is valid for $0<t<\kappa\Phi(\overline{R})$ with some $\kappa>0$ such that $(1/(2\phi^{-1}(t^{-1})),1/\phi^{-1}(\rho t^{-1}))\subset (0,\psi_c(\overline{R}))$ (where $\phi^{-1}(t^{-1})\asymp 1/\psi_c(\Phi^{-1}(c_\bullet t))$) for such $t$. Then, by the semigroup property, $\mathrm{LLE}(\Phi)$ holds for $0<t<\Phi(\overline{R})$. The implication $\mathrm{LLE}(\Phi)+\mathrm{UJS}\Rightarrow\mathrm{ABB}(\Phi)+\mathrm{PI}(\Phi)$ follows from Proposition \ref{PHI}. Since both $\mathcal{E}_j$ and $\mathcal{E}_{j,0}$ satisfy $\mathrm{J}(\psi_j)$, they are comparable. Thus $\mathrm{UJS}$, $\mathrm{ABB}(\Phi)$ and $\mathrm{PI}(\Phi)$ hold for $\mathcal{E}_j$, while $\mathrm{J}_\le(\Phi)$ follows from (\ref{comp2}). Again by Proposition \ref{PHI}, $\mathcal{E}_j$ satisfies $\mathrm{NLE}(\Phi)$. Based on this result, the proof to
$$p_t(x,y)\ge\frac{C^{-1}t}{V(x,d(x,y))\psi_j(d(x,y))}\mbox{ for all }d(x,y)>\eta\Phi^{-1}(t)$$
is nothing different from \cite[Part (ii) in the proof of Proposition 5.4]{ckw0}. Therefore, it suffices to prove that there exist constants $C,c,\kappa>0$ such that
\begin{equation}\label{S-left}
p_t(x,y)\ge\frac{C^{-1}}{V(x,\Phi^{-1}(t))}\exp\left(-c\sup\limits_{\sigma>0}\left[\frac{\rho(d(x,y))}{\rho(\sigma)}-\frac{\kappa t}{\Phi(\sigma)}\right]\right)\mbox{ if }d(x,y)>\eta\Phi^{-1}(t),
\end{equation}
and then apply Lemma \ref{3.8} to obtain an analogue without $\kappa$.

Actually, by \cite[Lemma 6.3]{gt}, for all $0<\varepsilon<d(x,y)$, there exists $x_1,\cdots,x_{N-1}\in M$, where 
$$d(x,x_1)<\varepsilon,\quad d(x_i,x_{i+1})<\varepsilon,\quad d(x_{N-1},y)<\varepsilon,\quad\mbox{and}\quad \frac{d_\varepsilon(x,y)}{\varepsilon}\le N<10\frac{d_\varepsilon(x,y)}{\varepsilon}.$$

Since $\beta_1^{(r_*)}(\Phi\circ\rho^{-1})>1$, then
$$\varlimsup\limits_{\varepsilon\to 0}\frac{\Phi(\varepsilon)}{\varepsilon}d_\varepsilon(x,y)\le\varlimsup\limits_{\varepsilon\to 0}\Phi(\varepsilon)\cdot C\frac{\rho(d(x,y))}{\rho(\varepsilon)}\xlongequal{\varepsilon=\rho(\delta)}C\rho(d(x,y))\varlimsup\limits_{\delta\to 0}\frac{\Phi(\rho^{-1}(\delta))}{\delta}\le C'\rho(d(x,y))\varlimsup\limits_{\delta\to 0}\delta^{\beta_1^{(r_*)}(\Phi\circ\rho^{-1})-1}=0.$$
Therefore, parallel to the proof of (\ref{opt-n_eps}), if $\varepsilon=\varepsilon'(t,x,y)$ is taken so that $\frac{\Phi(\varepsilon)}{\varepsilon}d_\varepsilon(x,y)=t\le\Phi(r_*)$, then
\begin{equation}\label{opt}
C_7^{-1}\frac{d_\varepsilon(x,y)}{\varepsilon}\le\sup\limits_{r>0}\left(\frac{\rho(d(x,y))}{\rho(r)}-\frac{t}{\Phi(r)}\right)\le C_7\frac{d_\varepsilon(x,y)}{\varepsilon}.
\end{equation}

Now taking $\varepsilon=\varepsilon'(\delta t,x,y)$, where
$$\delta:=\frac{1}{10C(\Phi)}\left(\frac{\eta\wedge 1}{3}\right)^{\beta_2(\Phi)}$$
so that $3\varepsilon\le\eta\Phi^{-1}(t/N)$, then for all $y_i\in B(x_i,\varepsilon)$, we have
$$d(y_i,y_{i+1})\le d(y_i,x_i)+d(x_i,x_{i+1})+d(x_{i+1},y_{i+1})<3\varepsilon\le\eta\Phi^{-1}(t/N),$$
and similarly $d(x,y_1)$, $d(y_{N-1},y)<\eta\Phi^{-1}(t/N)$. Thus by $\mathrm{NLE}(\phi)$, for all $t<\Phi(r_*\wedge\overline{R})$,
\begin{align*}
p_t(x,y)=&\int_{M^{N-1}}p_{t/N}(x,y_1)p_{t/N}(y_1,y_2)\cdots p_{t/N}(y_{N-1},y)d(\mu\otimes\cdots\otimes\mu)(y_1,y_2,\cdots,y_{N-1})\\
\ge&\int_{B(x_1,\varepsilon)}\cdots\int_{B(x_{N-1},\varepsilon)}p_{t/N}(x,y_1)p_{t/N}(y_1,y_2)\cdots p_{t/N}(y_{N-1},y)d\mu(y_{N-1})\cdots d\mu(y_1)\\
\ge&\int_{B(x_1,\varepsilon)}\cdots\int_{B(x_{N-1},\varepsilon)}\frac{C^{-1}}{V(x,\Phi^{-1}(t/N))}\prod\limits_{i=1}^{N-1}\frac{C^{-1}}{V(x_i,\Phi^{-1}(t/N))}d\mu(y_{N-1})\cdots d\mu(y_1)\\
=&\frac{C^{-1}}{V(x,\Phi^{-1}(t/N))}\prod\limits_{i=1}^{N-1}\left(\frac{C^{-1}}{V(x_i,\Phi^{-1}(t/N))}\cdot V(x_i,\varepsilon)\right)\ge\frac{C^{-1}}{V(x,\Phi^{-1}(t/N))}\prod\limits_{i=1}^{N-1}C^{-1}\\
\ge&\frac{1}{V(x,\Phi^{-1}(t))}\cdot\exp\{-\log(C)\cdot N\}\ge\frac{1}{V(x,\Phi^{-1}(t))}\exp\left\{-10\log(C)\frac{d_\varepsilon(x,y)}{\varepsilon}\right\}\\
\ge&\frac{1}{V(x,\Phi^{-1}(t))}\exp\left\{-10C_7\log(C)\sup\limits_{r>0}\left(\frac{\rho(d(x,y))}{\rho(r)}-\frac{\delta t}{\Phi(r)}\right)\right\},
\end{align*}
which is exactly (\ref{S-left}) with $C=1$, $c=10C_7\log C$ and $\kappa=\delta$. Hence the lower bound in (\ref{S+HK}) is proved.

(2) Now we turn to the upper bound. Clearly by (\ref{comp2}), there exists $C_j>0$ such that $\Phi\circ F_{\rho,\varepsilon}(r)\le C_j\psi_j\circ F_{\rho,\varepsilon}(r)$ for all $\varepsilon>0$ and $r>0$. Take $\hat{\Phi}(r):=K\Phi(r)$ and $\hat{\psi}_j(r):=C_jK\psi_j(r)$, with $K\ge 1$ to be determined. Then, $\mathrm{J}(\hat{\psi}_j)$ holds, implying  $\mathrm{J}_\le(\hat{\Phi})$. Further, we cover \cite[Theorem 4.6]{bkkl} with a bounded version: same as in Part (1), $\mathrm{UJS}$, $\mathrm{ABB}(\hat{\Phi})$ and $\mathrm{PI}(\hat{\Phi})$ holds; again by Proposition \ref{PHI}, $\mathrm{DUE}(\hat{\Phi})$ and $\mathrm{E}(\hat{\Phi})$ follow for $\mathcal{E}_j$.

For each $\varepsilon\in I:=(0,\infty)$, let $(M^{(\varepsilon)},d^{(\varepsilon)},\mu^{(\varepsilon)},\mathcal{E}^{(\varepsilon)},\mathcal{F}^{(\varepsilon)})=(M,d_\varepsilon,\mu,\mathcal{E}_j,\mathcal{F}_j)$, where conditions $\mathrm{VD}$, $\mathrm{RVD}$, $\mathrm{DUE}(\hat{\Phi}\circ F_{\rho,\varepsilon})$, $\mathrm{E}(\hat{\Phi}\circ F_{\rho,\varepsilon})$ and $\mathrm{J}(\hat{\psi}_j\circ F_{\rho,\varepsilon})$ all hold uniformly on $\varepsilon>0$  by Remark \ref{release-T} and Lemma \ref{key}, and $\hat{\Phi}(r)\le\hat{\psi}_j(r)$. Further, parallel to the proof of Corollary \ref{beta-uni'}, we see $\beta_1^{(r_{*,\varepsilon})}(\hat{\Phi}\circ F_{\rho,\varepsilon})\ge\beta_1^{(r_*)}(\hat{\Phi}\circ\rho^{-1})>1$, where $r_{*,\varepsilon}=F_{\rho,\varepsilon}^{-1}(r_*)$.

Consequently, by Proposition \ref{found+}, uniformly on $\varepsilon>0$, we have for all $x,y\in M$, $\varepsilon>0$ and $0<t<\hat{\Phi}(r_*)$,
\begin{align}
p_t(x,y)\le&\frac{C}{V_\varepsilon(x,(\Phi\circ F_{\rho,\varepsilon})^{-1}(t))}\wedge\left\{\frac{Ct}{V_\varepsilon(x,d_\varepsilon(x,y))(\psi_j\circ F_{\rho,\varepsilon})(d_\varepsilon(x,y))}+\right.\notag\\
&\quad\quad\quad\quad\quad\quad\quad\quad\quad\left.+\frac{C}{V_\varepsilon(x,(\Phi\circ F_{\rho,\varepsilon})^{-1}(t))}\exp\left(-c\sup\limits_{\sigma>0}\left[\frac{d_\varepsilon(x,y)}{\sigma}-\frac{t}{\Phi\circ F_{\rho,\varepsilon}(\sigma)}\right]\right)\right\}\notag\\
=&\frac{C}{V_\varepsilon(x,F_{\rho,\varepsilon}^{-1}(\Phi^{-1}(t)))}\wedge\left\{\frac{Ct}{V_\varepsilon(x,d_\varepsilon(x,y))\psi_j(F_{\rho,\varepsilon}(d_\varepsilon(x,y)))}+\right.\notag\\
&\quad\quad\quad\quad\quad\quad\quad\quad\quad\left.+\frac{C}{V_\varepsilon(x,F_{\rho,\varepsilon}^{-1}(\Phi^{-1}(t)))}\exp\left(-c\sup\limits_{\sigma>0}\left[\frac{d_\varepsilon(x,y)}{\sigma}-\frac{t}{\Phi(F_{\rho,\varepsilon}(\sigma))}\right]\right)\right\}\notag\\
\le&\frac{C}{V(x,C^{-1}\Phi^{-1}(t))}\wedge\left\{\frac{Ct}{V(x,C^{-1}d(x,y))\psi_j(Cd(x,y))}+\right.\notag\\
&\quad\quad\quad\quad\quad\quad\quad\quad\quad\left.+\frac{C}{V(x,C^{-1}\Phi^{-1}(t))}\exp\left(-c\sup\limits_{\sigma>0}\left[\frac{d_\varepsilon(x,y)}{\sigma}-\frac{t}{\Phi(F_{\rho,\varepsilon}(\sigma))}\right]\right)\right\}\notag\\
\le&\frac{C'}{V(x,\Phi^{-1}(t))}\wedge\left\{\frac{C't}{V(x,d(x,y))\psi_j(d(x,y))}+\frac{C'}{V(x,\Phi^{-1}(t))}\exp\left(-c\left[\frac{d_\varepsilon(x,y)}{\varepsilon}-\frac{t}{\Phi(\varepsilon)}\right]\right)\right\}.\label{uhk-2}
\end{align}
Parallel to Part (1) in the proof of Lemma \ref{key}, we see $\Phi(r_*\wedge\overline{R})\le\hat{C}\Phi\circ F_{\rho,\varepsilon}(r_{*,\varepsilon}\wedge\overline{R}_\varepsilon)$ independently of $\varepsilon>0$. Therefore, by taking $K=\hat{C}$ we have (\ref{uhk-2}) holds for all $0<t<\Phi(r_*\wedge\overline{R})$.

Taking the infimum of the right side of (\ref{uhk-2}) over $\varepsilon\in(0,d(x,y))$, we see
\begin{align*}
p_t(x,y)\le&\frac{C'}{V(x,\Phi^{-1}(t))}\wedge\left\{\frac{C't}{V(x,d(x,y))\psi_j(d(x,y))}+\frac{C'}{V(x,\Phi^{-1}(t))}\exp\left(-c\sup\limits_{0<\varepsilon<d(x,y)}\left[\frac{d_\varepsilon(x,y)}{\varepsilon}-\frac{t}{\Phi(\varepsilon)}\right]\right)\right\}\\
\le&\frac{C'}{V(x,\Phi^{-1}(t))}\wedge\left\{\frac{C't}{V(x,d(x,y))\psi_j(d(x,y))}+\frac{C'}{V(x,\Phi^{-1}(t))}\exp\left(-c'\sup\limits_{0<\varepsilon<d(x,y)}\left[\frac{\rho(d(x,y))}{\rho(\varepsilon)}-\frac{t}{\Phi(\varepsilon)}\right]\right)\right\}\\
\le&\frac{C'}{V(x,\Phi^{-1}(t))}\wedge\left\{\frac{C't}{V(x,d(x,y))\psi_j(d(x,y))}+\frac{C''}{V(x,\Phi^{-1}(t))}\exp\left(-c\sup\limits_{\varepsilon>0}\left[\frac{\rho(d(x,y))}{\rho(\varepsilon)}-\frac{t}{\Phi(\varepsilon)}\right]\right)\right\},
\end{align*}
where $c'=c/C(d,\rho)$ and $C''=C'e^{c'}$ for the third line since
$$\sup\limits_{\varepsilon\ge d(x,y)}\left[\frac{\rho(d(x,y))}{\rho(\varepsilon)}-\frac{t}{\Phi(\varepsilon)}\right]\ge\frac{\rho(d(x,y))}{\rho(d(x,y)\vee\Phi^{-1}(t))}-\frac{t}{\Phi(d(x,y)\vee\Phi^{-1}(t))}\ge -1.$$
Hence the upper bound in (\ref{S+HK}) is proved.
\end{proof}

Here are some remarks:

In the unbounded case, the subordinated diffusion $Y_t=X_{S_t}$, where the subordinator $S$ corresponds to L\'evy measure
$$\nu(dt)=\frac{dt}{t\psi_j\circ\psi_c^{-1}(t)},$$
satisfies $\mathrm{J}(\psi_j)$ and $\mathrm{PHI}(\Phi)$,  provided that $\int_0^1t\nu(dt)<\infty$ (cf. \cite[Theorem 2.3 and Corollary 2.5]{lm}). Clearly its heat kernel exists and by applying the known distribution law of $S_t$ (see \cite{mi}), it satisfies
\begin{equation}\label{ghk-both}
p_t(x,y)\asymp\frac{1}{V(x,\Phi^{-1}(t))}\wedge\left\{\frac{t}{V(x,d(x,y))\psi_j(d(x,y))}+\frac{1}{V(x,\Phi^{-1}(t))}\exp\left(-c_\bullet\sup\limits_{\sigma>0}\left[\frac{\rho(d(x,y))}{\rho(\sigma)}-\frac{\psi_c(\Phi^{-1}(t))}{\psi_c(\sigma)}\right]\right)\right\}.
\end{equation}
Clearly by Theorem \ref{thm2}, (\ref{ghk-both}) is equivalent to (\ref{S+HK}) if $\beta_1^{(r_*)}(\Phi\circ\rho^{-1})>1$, and thus holds for any jump kernel satisfying $\mathrm{J}(\psi_j)$ with the same $\psi_j$. Actually, even without the assumption $\beta_1^{(r_*)}(\Phi\circ\rho^{-1})>1$, by applying \cite[Theorem 2.19]{bkkl} instead of \cite[Corollary 2.15]{bkkl} in the proof above, we can directly obtain the upper bound in (\ref{ghk-both}), which, together with the lower bound in G$^+$HK$(\Phi,\psi_j)$, is called GHK$(\Phi,\psi_j)$. However, we can not obtain the lower bound of (\ref{ghk-both}) as in Part (1), since the adjacent points $y_i,y_{i+1}$ in an $\varepsilon$-chain, where $\varepsilon$ is nearly optimal for the supremum in this formula, are not near-diagonal (due to (\ref{comp3})) at time $t/N$.

For a general metric measure space that admits a diffusion satisfying $\mathrm{PHI}(\psi_c)$, clearly there exist scale functions $\rho_1,\rho_2$ (for example, $\rho_1(r)=r$ and $\rho_2(r)=\sqrt{\psi_c(r)}$) such that for all $x,y\in M$ and $0<\varepsilon<d(x,y)$,
\begin{equation}
C^{-1}\frac{\rho_1(d(x,y))}{\rho_1(\varepsilon)}\le\frac{d_\varepsilon(x,y)}{\varepsilon}\le C\frac{\rho_2(d(x,y))}{\rho_2(\varepsilon)}.
\end{equation}
This is a weaker chain condition. In this case, the method we use to prove Theorem \ref{thm2} can not yield the expected result
$$p_t(x,y)\asymp\frac{1}{V(x,\Phi^{-1}(t))}\wedge\left\{\frac{t}{V(x,d(x,y))\psi_j(d(x,y))}+\frac{1}{V(x,\Phi^{-1}(t))}\exp\left(-c_\bullet\sup\limits_{\varepsilon>0}\left[\frac{d_\varepsilon(x,y)}{\varepsilon}-\frac{t}{\Phi(\varepsilon)}\right]\right)\right\},$$
since we can only obtain $\mathrm{DUE}(\Phi\circ F_{\rho_1,\varepsilon})$, $\mathrm{E}_\le(\Phi\circ F_{\rho_1,\varepsilon})+\mathrm{E}_\ge(\Phi\circ F_{\rho_2,\varepsilon})$ and $\mathrm{J}_\le(\psi\circ F_{\rho_1,\varepsilon})+\mathrm{J}_\ge(\psi\circ F_{\rho_2,\varepsilon})$ under the metric $d_\varepsilon$. Actually, it is still an open question whether there exists a sharp stable heat kernel estimate whose two sides match in the general case.

\section{Examples for the $\rho$-chain condition}
Let us begin with a trivial example. Let $(M,d)$ be a metric space satisfying the chain condition. Then, $(M,d^\gamma)$, where $d^\gamma(x,y):=(d(x,y))^\gamma$ with an arbitrary $0<\gamma<1$, satisfies the $\rho_\gamma$-chain condition (where $\rho_\gamma:r\mapsto r^{1/\gamma}$). Further, if there exists a diffusion satisfying a two-sided sub-Gaussian estimate
$$p_t(x,y)\asymp\frac{1}{V(x,\psi_c^{-1}(t))}\exp\left\{-c_\bullet\sup\limits_{\sigma>0}\left(\frac{d(x,y)}{\sigma}-\frac{t}{\psi_c(\sigma)}\right)\right\},$$
then by letting $\delta=\rho_\gamma^{-1}(\varepsilon)=\varepsilon^\gamma$ we have
\begin{align*}
p_t(x,y)\asymp&\frac{1}{V_{d^\gamma}(x,(\psi_c\circ\rho_\gamma)^{-1}(t))}\exp\left\{-c_\bullet\sup\limits_{\sigma>0}\left(\frac{(d^\gamma(x,y))^{1/\gamma}}{\sigma}-\frac{t}{\psi_c(\sigma)}\right)\right\}\\
=&\frac{1}{V_{d^\gamma}(x,(\psi_c\circ\rho_\gamma)^{-1}(t))}\exp\left\{-c_\bullet\sup\limits_{\delta>0}\left(\frac{\rho_\gamma(d^\gamma(x,y))}{\rho_\gamma(\delta)}-\frac{t}{\psi_c\circ\rho_\gamma(\delta)}\right)\right\},
\end{align*}
which is exactly the estimate (\ref{rho-gauss}) with the scale $\psi_c\circ\rho_\gamma$. In the same way we can obtain the estimate (\ref{S+HK}) with the scales $\psi_j\circ\rho_\gamma$ and $\Phi\circ\rho_\gamma$ from (\ref{shk}) directly.

Now we bring up with a non-trivial example. Consider the Sierpi\'nski gasket (for short, SG). Let $p_1=\frac{1}{2}(1+\sqrt{-3}),p_2=0,p_3=1\in\mathbb{C}$, and $F_1,F_2,F_3$ be the contractions on $\mathbb{C}$ with $F_i(x):=\frac{1}{2}(x+p_i)$. SG is defined as the unique compact set $K\subset\mathbb{C}$ such that $K=F_1(K)\cup F_2(K)\cup F_3(K)$, and its essential boundary is $V_0:=\{p_1,p_2,p_3\}$. Any resistance form on $V_0$ can be represented by
$$E(u)=\min\limits_{y\in\mathbb{R}}\left\{a_1(u(p_1)-y)^2+a_2(u(p_2)-y)^2+a_3(u(p_3)-y)^2\right\}.$$
Kigami's renormalization equation reads
$$E(u)=\min\bigg\{\sum_{i=1}^3w_iE(v\circ F_i):\ v\mbox{ is a function on }V_1=F_1(V_0)\cup F_2(V_0)\cup F_3(V_0),\mbox{ and }v|_{V_0}=u\bigg\},$$
where $w_i>0$ are called the self-similar weights. The solution is called regular, if every $w_i>1$. In this case, let $\Lambda_r=\{\omega=i_1\cdots i_n:i_k\in\{1,2,3\}\mbox{ for all }1\le k\le n\mbox{ and }w_{i_1}^{-1}\cdots w_{i_n}^{-1}<r\le w_{i_1}^{-1}\cdots w_{i_{n-1}}^{-1}\}$ be the partition of $K$ with level $r\in(0,1)$. Let $F_\omega=F_{i_1}\circ\cdots\circ F_{i_n}$ and $w_\omega=w_{i_1}\cdots w_{i_n}$, and then the limit form
$$\mathcal{E}(u)=\lim\limits_{r\to 0}\sum\limits_{\omega\in\Lambda_r}w_\omega E(u\circ F_\omega)\mbox{ for all }u\in C(K)$$
expands to a regular strongly local conservative Dirichlet form which admits a resistance metric $R$ (depending on $t$) on $K$ such that the topology induced by $R$ is equal to the original topology on $K$. See \cite[Theorem 2.10]{hk}.

We consider only the resistance forms and weights of type $a_1=\lambda^{-1}>0$, $a_2=a_3=1$, $w_1=\sigma^{-1}$ and $w_2=w_3=\tau^{-1}$, with $0<\sigma,\tau<1$. Then, by the $\Delta$-Y transform, $(E,\{w_i\}_{i=1}^3)$ is a solution to the renormalization problem if and only if $\lambda,\sigma,\tau$ satisfy the following system of equations:
$$\lambda\sigma+\frac{(\sigma+\lambda\tau)^2}{2(\sigma+\tau+\lambda\tau)}=\lambda,\ \mbox{ and }\ \tau+\frac{\tau(\sigma+\lambda\tau)}{\sigma+\tau+\lambda\tau}=1,$$
or equivalently,
\begin{align}
\lambda\tau(1+\tau-2\sigma)-(1-\tau)\sigma&=0,\label{lambda-and-t}\\
f_\tau(\sigma):=2(2\tau-1)\sigma^2+2(\tau^2-3\tau+1)\sigma+\tau(1-\tau^2)&=0.\label{s-and-t}
\end{align}
For any $0<\tau<1$, since $f_\tau(1-\tau)=\tau(1-\tau)^2>0$ and $f_\tau(1)=-\tau(1-\tau)^2<0$, there exists a unique $\sigma=\sigma(\tau)\in (1-\tau,1)$ satisfying (\ref{s-and-t}). Additionally, note that $f_\tau(\frac{1}{2}(1+\tau))=(\tau+1)(\tau-1)(\tau-\frac{1}{2})$. Therefore, if $0<\tau\le\frac{1}{2}$, then $\sigma(\tau)\ge\frac{1}{2}(1+\tau)$, and thus the $\lambda$ satisfying (\ref{lambda-and-t}) is negative, not admissible for a resistance form; while if $\frac{1}{2}<\tau<1$, $\sigma(\tau)\in(1-\tau,\frac{1}{2}(1+\tau))$, and a corresponding renormalized resistance form follows. Note that in this case, $\sigma+\tau>1$.

By \cite[Lemma 5.3]{hk}, the near-diagonal lower estimate holds:
$$p_t(x,y)\ge\frac{C^{-1}}{t^{\alpha_*/\beta_*}},\mbox{ for all }0<t\le 1\mbox{ and }x,y\in K\mbox{ such that }t>R(x,y)^{\beta_*},$$
where $\alpha_*$ satisfies the Moran equation $\sigma^{\alpha_*}+\tau^{\alpha_*}+\tau^{\alpha_*}=1$, and $\beta_*=\alpha_*+1$. Actually the sub-Gaussian upper estimate also holds, see the paragraph after \cite[Corollary 5.3]{liu}.

Following \cite[Lemma 3.7 and Corollary 5.3]{liu} (or, in a discrete form, \cite[Lemma 3.2]{hk}), we see for all $0<r<1$ that $R(x,y)\le Cr$ for every two points $x,y\in K$ in a same $r$-cell $K_\omega$ (that is, $\omega\in\Lambda_r$), and $R(x,y)\ge C^{-1}r$ for every two points $x,y\in K$ that do not lie on adjacent $r$-cells (that is, there do not exist $\omega,\omega'\in\Lambda_r$ such that $x\in K_\omega,y\in K_{\omega'}$ and $K_\omega\cap K_{\omega'}\ne\emptyset$).

\begin{lemma}\label{example-lem}
The metric space $(K,R)$ does not satisfy the chain condition. For $\frac{3}{5}<\tau<1$, there does not exist a scale function $\rho$ such that $(K,R)$ satisfies the $\rho$-chain condition; while for $\frac{1}{2}<\tau\le\frac{3}{5}$, there exists a scale function $\rho:r\mapsto r^\gamma$ with $1<\gamma<\beta_*$, such that $(K,R)$ satisfies the $\rho$-chain condition.
\end{lemma}

\begin{proof}
(1) By \cite[Lemma 6.2]{hk}, letting $\gamma_1,\gamma_2$ be the solutions to
\begin{equation}\label{side-moran}
\sigma^{\gamma_1}+\tau^{\gamma_1}=1\ \mbox{ and }\ \tau^{\gamma_2}+\tau^{\gamma_2}=1,
\end{equation}
we have: if $\gamma_1\le\gamma_2$, then $R_r(p_i,p_j)\asymp r\cdot r^{-\gamma_1}$ for any $\{i,j\}\subset\{1,2,3\}$; while if $\gamma_1>\gamma_2$, then $R_r(p_2,p_3)\asymp r\cdot r^{-\gamma_2}$ but $R_r(p_1,p_i)\asymp r\cdot r^{-\gamma_1}$ for $i=2$ or $3$.

Since $0<\sigma<1,\frac{1}{2}<\tau<1$ and $\sigma+\tau>1$, then $\gamma_1,\gamma_2>1$. Therefore, $R_r(p_i,p_j)\to\infty$ as $r\to 0$, which implies that $(K,R)$ does not satisfy the chain condition.

(2) It is clear by (\ref{lambda-and-t}), (\ref{s-and-t}) and (\ref{side-moran}) that
$$\gamma_1\le\gamma_2\Leftrightarrow \sigma\ge\tau\Leftrightarrow\lambda\ge 1\Leftrightarrow\frac{1}{2}<\tau\le\frac{3}{5}.$$
Therefore, for $\frac{3}{5}<\tau<1$, $R_r(p_2,p_3)$ and $R_r(p_1,p_2)$ grows in different races, which implies that there is no scale $\rho$ such that $(K,R)$ satisfies the $\rho$-chain condition.

(3) For $\frac{1}{2}<\tau\le\frac{3}{5}$, $R_r(p_i,p_j)\asymp r^{1-\gamma_1}$. Consequently, for two arbitrary points $x,y\in K$, let $r_0$ be the smallest $r>0$ such that $x,y$ are in adjacent $r$-cells, say $x\in K_\omega,y\in K_{\omega'}$ and $\{q\}=K_\omega\cap K_{\omega'}$. Then $r_0\asymp R(x,y)$. Without loss of generality, we assume $x$ is not in a same $(r_0\sqrt{\sigma\vee\tau})$-cell (that is, $q=F_\omega(p_i)$ and $x\in K_{\omega j}$ with some $j\ne i$). Note that any $r$-chain from $x$ to $y$ must pass $B_R(q,r)$, which implies that
$$R_r(x,q)+R_r(q,y)\ge R_r(x,y)\ge\inf\limits_{\xi\in B_R(q,r)}R_r(x,\xi)+\inf\limits_{\xi\in B_R(q,r)}R_r(\xi,y)\ge R_r(x,q)-r+0>\frac{1}{2}R_r(x,q),$$
if $r<(2C)^{-1}(\sigma\wedge\tau)r_0$ (since $r<\frac{1}{2}R(x,q)\le\frac{1}{2}R_r(x,q)$ in this case). Let $q$ take the place of $y$ and repeat this step, we have
$$R_r(x,q_1)+R_r(q_1,q)\ge R_r(x,q)\ge\frac{1}{2}R_r(q_1,q)$$
with $r<(2C)^{-1}(\sigma\wedge\tau)^2r_0$, where $\{q_1\}=K_{\omega i}\cap K_{\omega j}$. Therefore,
$$R_r(x,q)\ge\frac{1}{2}R_r(q_1,q)=\frac{1}{2}R_r(F_{\omega i}(p_k),F_{\omega i}(p_i))=\frac{1}{2}w_{\omega i}^{-1}R_{w_{\omega i}r}(p_i,p_j)\asymp r_0(r_0/r)^{\gamma_1-1}=r(r_0/r)^{\gamma_1},$$
and by letting $q_1$ take the place of $q$ and get $q_2$ form the same procedure and so on, until some $q_n$ lies on the same $C^{-1}r$-cell with $x$ (so that $R_r(x,q_n)=R(x,q_n)\le C\cdots C^{-1}r=r$), we have
$$R_r(x,q)\le\sum\limits_{i=0}^{n-1}R_r(q_i,q_{i+1})+r\lesssim\sum\limits_{i=0}^{n-1}r((\sigma\vee\tau)^ir_0/r)^{\gamma_1}=r(r_0/r)^{\gamma_1}\sum\limits_{i=0}^{n-1}(\sigma\vee\tau)^{i\gamma_1}\lesssim r(r_0/r)^{\gamma_1}.$$
The same is true for $R_r(q,y)$. Therefore,
$$R_r(x,y)\asymp r(r_0/r)^{\gamma_1}\asymp r\left(\frac{R(x,y)}{r}\right)^{\gamma_1},$$
that is, when $\frac{1}{2}<\tau\le\frac{3}{5}$, $(K,R)$ satisfies the $\rho$-chain condition, where $\rho(r)=r^{\gamma_1}$. Further, since
$$\sigma^{\alpha_*}+\tau^{\alpha_*}<\sigma^{\alpha_*}+\tau^{\alpha_*}+\tau^{\alpha_*}=1,$$
then by (\ref{side-moran}), $\gamma_1<\alpha_*$, and thus $\gamma_1<\alpha_*+1=\beta_*$. The lemma follows by taking $\gamma=\gamma_1$.
\end{proof}

After proper dilation, we assume that $diam(K,R)=1$. By Theorem \ref{thm1}:

\begin{corollary}
For any $\frac{1}{2}<\tau\le\frac{3}{5}$, the corresponding self-similar Dirichlet form admits the following two-sided heat kernel estimate:
$$p_t(x,y)\asymp\frac{1}{t^{\alpha_*/\beta_*}}\exp\left\{-c_\bullet\left(\frac{R(x,y)^{\beta_*}}{t}\right)^{\frac{\gamma}{\beta_*-\gamma}}\right\},$$
for all $0<t<1$ and $x,y\in K$.\qed
\end{corollary}

This is the same estimate as in \cite[Part (1) of Theorem 6.3]{hk}.

Now we turn to jump-type forms. As follows is a typical example where $\Phi$ is not comparable with $\psi_j$, and the effect of the $\rho$-chain condition can be seen by comparing this estimate with \cite[Example 6.4]{bkkl0}:

\begin{corollary}\label{jump-sg}
For any $\frac{1}{2}<\tau\le\frac{3}{5}$, let $J$ be a jump kernel satisfying $\mathrm{J}(\psi_j)$, where $\psi_j(r)=r^{\beta_*}\log^2(1/(r\wedge e^{-2/\beta_*}))$, then $J$ defines a regular Dirichlet form satisfying the following heat kernel estimate:
\begin{align*}
p_t(x,y)\asymp\left(\frac{\log(t^{-1}\vee e^{2/\beta_*})}{t}\right)^{\frac{\alpha_*}{\beta_*}}\wedge&\left(\frac{t}{R(x,y)^{\alpha_*+\beta_*}\log^2(R(x,y)^{-1}\vee e^{2/\beta_*})}+\right.\\
&+\left.\left(\frac{\log(t^{-1}\vee e^{2/\beta_*})}{t}\right)^{\frac{\alpha_*}{\beta_*}}\exp\left\{-c_\bullet\left(\frac{R(x,y)^{\beta_*}}{t}\log\left(\frac{R(x,y)^\gamma}{t}\vee e^{2/\beta_*}\right)\right)^{\frac{\gamma}{\beta_*-\gamma}}\right\}\right),
\end{align*}
for all $0<t<1$ and $x,y\in K$.
\end{corollary}
Here the parameter $e^{-2/\beta_*}$ is taken so that $\psi_j$ is monotonic.
\begin{proof}
By elementary calculus we have for each $0<r\le e^{-2/\beta_*}$,
$$\Phi(r)=\frac{\psi_c(r)}{\int_0^r\frac{\psi_c(s)ds}{s\psi_j(s)}}=\frac{r^{\beta_*}}{\int_0^r\frac{s^{\beta_*}ds}{s\cdot s^{\beta_*}\log^2(1/s)}}\xlongequal{s=e^t}\frac{r^{\beta_*}}{\int_{-\infty}^{\log r}\frac{dt}{t^2}}=\frac{r^{\beta_*}}{(-\log r)^{-1}}=r^{\beta_*}\log(1/r),$$
and clearly $\Phi(r)\asymp r^{\beta_*}$ for $r>e^{-2/\beta_*}$. Without loss of generality, we take
$$\Phi(r)=r^{\beta_*}\log(1/(r\wedge e^{-2/\beta_*})).$$
This way, if $r=\Phi^{-1}(t)$, then $r^{\beta_*}\log(1/(r\wedge e^{-2/\beta_*}))=t$. If $0<r\le e^{-2/\beta_*}$, then 
\begin{equation}\label{theta-1}
r^{\beta_*}\log(1/r)=t.
\end{equation}
Therefore,
$$\beta_*\log(r^{-1})-\log\log(r^{-1})=\log(t^{-1}).$$
Since $\beta_*\ge 2>2/e$, then $0<\log\log(r^{-1})\le \frac{\beta_*}{2}\log(r^{-1})$ for all $0<r\le e^{-2/\beta_*}$. Therefore, we have
\begin{equation}\label{theta-2}
\frac{\beta_*}{2}\log(r^{-1})\le\log(t^{-1})\le\beta_*\log(r^{-1}).
\end{equation}
Combining (\ref{theta-1}) with (\ref{theta-2}), we obtain for all $0<r\le e^{-2/\beta_*}$ that
\begin{equation}\label{theta-3}
\left(\frac{\beta_*t}{2\log(t^{-1})}\right)^{1/\beta_*}\le r=\Phi^{-1}(t)\le\left(\frac{\beta_*t}{\log(t^{-1})}\right)^{1/\beta_*}.
\end{equation}
In particular, if we take $K>1$ such that $K\log K\ge e^2\beta_*$, then for all $0<t\le K^{-1}$, we have
$$r\le\left(\frac{\beta_*t}{\log(t^{-1})}\right)^{1/\beta_*}\le\left(\frac{K^{-1}\beta_*}{\log K}\right)^{1/\beta_*}\le(e^{-2})^{1/\beta_*}=e^{-2/\beta_*}.$$
Therefore, (\ref{theta-3}) holds for all $0<t\le K^{-1}$. On the other hand, clearly $r=\Phi^{-1}(t)$ is almost a constant for $(K\beta_*)^{-1}<t\le 1$. Therefore, for all $0<t\le 1$,
$$\Phi^{-1}(t)\asymp\left(\frac{t\wedge K^{-1}}{\log((t\wedge K^{-1})^{-1})}\right)^{1/\beta_*}\asymp\left(\frac{t}{\log(t^{-1}\vee K)}\right)^{1/\beta_*}\asymp\left(\frac{t}{\log(t^{-1}\vee e^{2/\beta_*})}\right)^{1/\beta_*}.$$
Now let $\varepsilon=\varepsilon'(t;x,y)$ be the largest root of $\frac{\Phi(\varepsilon)}{\varepsilon}R_\varepsilon(x,y)=t$.
This way,
$$\varepsilon^{\beta_*-\gamma}\log(1/(\varepsilon\wedge e^{-2/\beta_*}))\asymp\frac{t}{R(x,y)^\gamma},$$
which implies if $t\le R(x,y)^\gamma$,
\begin{equation}\label{theta-4}
\varepsilon\asymp\left(\frac{t}{R(x,y)^\gamma}\bigg/\log\left(\frac{R(x,y)^\gamma}{t}\vee e^{2/\beta_*}\right)\right)^{1/(\beta_*-\gamma)}.
\end{equation}
If $t>R(x,y)^\gamma$, then obviously $\varepsilon>e^{-2/\beta_*}$, and thus we directly have $\varepsilon\asymp(t/R(x,y)^\gamma)^{1/(\beta_*-\gamma)}.$ Thus in each case (\ref{theta-4}) holds.

Further, by (\ref{opt}), with $\varepsilon$ defined as above,
$$\sup\limits_{\sigma>0}\left(\frac{\rho(R(x,y))}{\rho(\sigma)}-\frac{t}{\Phi(\sigma)}\right)\asymp\frac{R_\varepsilon(x,y)}{\varepsilon}\asymp\left(\frac{R(x,y)}{\varepsilon}\right)^\gamma.$$
Consequently, for all $x,y\in M$ and $0<t\le 1$,
\begin{align*}
p_t(x,y)\asymp&\frac{1}{V(x,\Phi^{-1}(t))}\wedge\left\{\frac{t}{V(x,R(x,y))\psi_j(R(x,y))}+\frac{1}{V(x,\Phi^{-1}(t))}\exp\left(-c_\bullet\sup\limits_{\sigma>0}\left[\frac{\rho(R(x,y))}{\rho(\sigma)}-\frac{t}{\Phi(\sigma)}\right]\right)\right\}.\\
\asymp&\left(\frac{\log(t^{-1}\vee e^{2/\beta_*})}{t}\right)^{\frac{\alpha_*}{\beta_*}}\wedge\left(\frac{t}{R(x,y)^{\alpha_*+\beta_*}\log^2(1/(R(x,y)\wedge e^{2/\beta_*}))}+\right.\\
&\quad\quad\quad\quad+\left.\left(\frac{\log(t^{-1}\vee e^{2/\beta_*})}{t}\right)^{\frac{\alpha_*}{\beta_*}}\exp\left\{-c'_\bullet\left(R(x,y)\left(\frac{t}{R(x,y)^\gamma}\bigg/\log\left(\frac{R(x,y)^\gamma}{t}\vee e^{2/\beta_*}\right)\right)^{-1/(\beta_*-\gamma)}\right)^\gamma\right\}\right)\\
\asymp&\left(\frac{\log(t^{-1}\vee e^{2/\beta_*})}{t}\right)^{\frac{\alpha_*}{\beta_*}}\wedge\left(\frac{t}{R(x,y)^{\alpha_*+\beta_*}\log^2(R(x,y)^{-1}\vee e^{2/\beta_*})}+\right.\\
&\quad\quad\quad\quad\quad\quad\quad\quad+\left.\left(\frac{\log(t^{-1}\vee e^{2/\beta_*})}{t}\right)^{\frac{\alpha_*}{\beta_*}}\exp\left\{-c'_\bullet\left(\frac{R(x,y)^{\beta_*}}{t}\log\left(\frac{R(x,y)^\gamma}{t}\vee e^{2/\beta_*}\right)\right)^{\frac{\gamma}{\beta_*-\gamma}}\right\}\right),
\end{align*}
which is the desired result.
\end{proof}

To end with, we remark on the case $\frac{3}{5}<\tau<1$ (that is, $\gamma_1>\gamma_2$): for any $x,y\in K$, take $r_0,\omega,\omega'$ as in Part (3) in the proof of Lemma \ref{example-lem}. Let $q=F_\omega(p_i)=F_{\omega'}(p_k)$, then $i\ne k$, say $i>k$. Let
$$L_x(q)=\{\xi\in K_\omega:\xi-q\in\mathbb{R}\},\quad L_y(q)=\{\xi'\in K_\omega:\xi'-q\in\mathbb{R}\}$$
(here we consider $K$ and $\mathbb{R}$ as subsets of $\mathbb{C}$), and denote
$$R^-(x,q):=\inf\{R(x,\xi):\xi\in L_x(q)\},\quad R^-(y,q):=\inf\{R(y,\xi):\xi\in L_y(q)\}.$$

Say $x\in K_{\omega j}$ and $y\in K_{\omega'l}$, then it is impossible if $j=i,l=k$. Therefore, the following four cases cover all possibilities of $i,j,k,l$:

(i) $\{i,j\}=\{2,3\}$ and $k=l=1$. Here
$$\frac{R_r(x,q)}{r}\asymp\inf\limits_{\xi\in L_x(q)}\left(\frac{R_r(x,\xi)}{r}+\frac{R_r(\xi,q)}{r}\right)\asymp 1\vee\left(\frac{R^-(x,q)}{r}\right)^{\gamma_1}+\left(\frac{r_0}{r}\right)^{\gamma_2},\quad\frac{R_r(q,y)}{r}\asymp\left(\frac{R(q,y)}{r}\right)^{\gamma_1}=\left(\frac{R^-(y,q)}{r}\right)^{\gamma_1};$$

(ii) $i=3,k=2$ and $j,l\in\{2,3\}$. Similar to (i), we obtain
$$1\vee\left(\frac{R^-(\xi,q)}{r}\right)^{\gamma_1}\lesssim\frac{R_r(\xi,q)}{r}\lesssim 1\vee\left(\frac{R^-(\xi,q)}{r}\right)^{\gamma_1}+\left(\frac{r_0}{r}\right)^{\gamma_2},\mbox{ where }\xi=x\mbox{ or }y,$$
and for at least one in $x$ and $y$, the second $\lesssim$ becomes $\asymp$;

(iii) $k=1$, $l\ne 1$ or $l=1$, $k\ne 1$. Similar to Step (2), we have
$$(R(x,y)/r)^{\gamma_1}\gtrsim\frac{R_r(x,y)}{r}\ge\frac{1}{2}\cdot\frac{R_r(q,y)}{r}\gtrsim(R(q,y)/r)^{\gamma_1}\gtrsim(R(x,y)/r)^{\gamma_1},\mbox{ and }R(x,y)\asymp R^-(y,q);$$

(iv) $j=1$. Similar to (iii), we have
$$(R(x,y)/r)^{\gamma_1}\gtrsim\frac{R_r(x,y)}{r}\ge\frac{1}{2}\cdot\frac{R_r(x,q)}{r}\gtrsim(R(x,q)/r)^{\gamma_1}\gtrsim(R(x,y)/r)^{\gamma_1},\mbox{ and }R(x,y)\asymp R^-(x,q).$$

In summary, in any case
\begin{align*}
\frac{R_r(x,y)}{r}\asymp&\left(\frac{R(x,y)}{r}\right)^{\gamma_2}\vee\left(\frac{R^-(x,q)}{r}\right)^{\gamma_1}\vee\left(\frac{R^-(y,q)}{r}\right)^{\gamma_1}\vee 1\\
\asymp&\begin{cases}
\left(\frac{R^-(x,q)\vee R^-(y,q)}{r}\right)^{\gamma_1},&\mbox{if }0<r\le r_\infty(x,y);\\
\left(\frac{R(x,y)}{r}\right)^{\gamma_2},&\mbox{if }r_\infty(x,y)<r<R(x,y),
\end{cases}\mbox{ where }r_\infty(x,y)=\frac{(R^-(x,q)\vee R^-(y,q))^{\frac{\gamma_1}{\gamma_1-\gamma_2}}}{R(x,y)^{\frac{\gamma_2}{\gamma_1-\gamma_2}}}.
\end{align*}
Consequently, the self-similar diffusion satisfies the following heat kernel estimate:
$$p^c_t(x,y)\asymp\frac{1}{t^{\alpha_*/\beta_*}}\exp\left\{-c_\bullet\left(\left(\frac{R^-(x,q)\vee R^-(y,q)}{t^{1/\beta_*}}\right)^{\frac{\beta_*\gamma_1}{\beta_*-\gamma_1}}\vee\left(\frac{R(x,y)}{t^{1/\beta_*}}\right)^{\frac{\beta_*\gamma_2}{\beta_*-\gamma_2}}\right)\right\},$$
which is a better estimate than \cite[Part (2) of Theorem 6.3]{hk}. Note that fixing $x,y\in K$, the second term is principal if and only if $t\ge[(R^-(x,q)\vee R^-(y,q))^{\gamma_1(\beta_*-\gamma_2)}/R(x,y)^{\gamma_2(\beta_*-\gamma_1)}]^{1/(\gamma_1-\gamma_2)}$.

Note also that the subordination techniques can not directly yield a heat kernel estimate for the subordinated process here, since the strongly local form only satisfy the sub-Gaussian upper bound for $0<t<1$ rather than for all $t>0$.

\section*{Acknowledgments}
The author would like to thank Professor Mathav Murugan for introducing the research topic, and Professor Alexander Grigor'yan for discussions on the examples.

\quad\\
Fakult\"at f\"ur Mathematik, Universit\"at Bielefeld, Postfach 100131, D-33501 Bielefeld, Germany\\
E-mail: gliu@math.uni-bielefeld.de


\begin{thebibliography}{ }
\bibitem{ab15}
S. Andres, M. T. Barlow. Energy inequalities for cutoff functions and some applications. J. Reine Angew. Math. 699 (2015) 183--215
\bibitem{bkkl0}
J. Bae, J. Kang, P. Kim, J. Lee. Heat kernel estimates for symmetric jump processes with mixed polynomial growths. Ann. Probab. 47 (2019) 2830--2868
\bibitem{bkkl}
J. Bae, J. Kang, P. Kim, J. Lee. Heat kernel estimates and their stabilities for symmetric jump processes with general mixed polynomial growths on metric measure spaces. To appear. %in: Proc. Lond. Math. Soc.
\bibitem{bbk}
M. T. Barlow, R. F. Bass, T. Kumagai. Stability of parabolic Harnack inequalities on metric measure spaces. J. Math. Soc. Japan 58 (2006) 485--519
\bibitem{bghh}
A. Bendikov, A. Grigor’yan, E. Hu, J. Hu. Heat kernels and non-local Dirichlet forms on ultrametric spaces. Ann. Sc.
Norm. Super. Pisa Cl. Sci. (5) 22(1) (2021) 399--461
\bibitem{ckw0}
Z.-Q. Chen, T. Kumagai, J. Wang. Stability of heat kernel estimates for symmetric jump processes on metric measure spaces. In press in Memoirs of the AMS. (2016)
\bibitem{ckw}
Z.-Q. Chen, T. Kumagai, J. Wang. Stability of parabolic Harnack inequalities for symmetric non-local Dirichlet forms. J. Eur. Math. Soc. 22 (2020) 3747--3803
\bibitem{fot}
M. Fukushima, Y. Oshima, M. Takeda. Dirichlet forms and symmetric Markov processes. 2nd version. Walter de Gruyter \& Co., Berlin (2011)
%\bibitem{gh08}
%A. Grigor'yan, J. Hu. Off-diagonal upper estimates for the heat kernel of the Dirichlet forms on metric spaces. Invent. Math. 174 (2008) 81--126
%\bibitem{ghh17}
%A. Grigor'yan, E. Hu, J. Hu. Lower estimates of heat kernels for non-local Dirichlet forms on metric measure spaces. J. Func. Anal. 272 (2017) 3311--3346
%\bibitem{ghh+}
%A. Grigor'yan, E. Hu, J. Hu. Tail estimates of heat kernels on doubling spaces. To appear.
\bibitem{ghl14}
A. Grigor'yan, J. Hu, K.-S. Lau. Estimates of heat kernels for non-local regular Dirichlet forms. Trans. Amer. Math. Soc. 366(12) (2014) 6397--6441
\bibitem{ghl15}
A. Grigor'yan, J. Hu, K.-S. Lau. Generalized capacity, Harnack inequality and heat kernels of Dirichlet forms on metric measure spaces. J. Math. Soc. Japan 67 (2015) 1485--1549
\bibitem{gt}
A. Grigor'yan, A. Telcs. Two-sided estimates of heat kernels on metric measure spaces. Ann. Probab. 40 (2012) 1212--1284
\bibitem{hk}
B. M. Hambly, T. Kumagai. Transition density estimates for diffusion processes on post critically finite self-similar fractals. Proc. Lond. Math. Soc. (3) 78 (1999) 431--458
\bibitem{hl}
J. Hu, G. Liu. Upper estimates of heat kernels for non-local Dirichlet forms on doubling spaces. Forum Math. 34 (2022) 225--277
\bibitem{liu}
G. Liu. Existence of self-similar Dirichlet forms on post-critically finite fractals in terms of their resistances. To appear.
\bibitem{lm}
G. Liu, M. Murugan. On the comparison between jump processes and subordinated diffusions. To appear.
\bibitem{mi}
A. Mimica. Heat kernel estimates for subordinate Brownian motions. Proc. Lond. Math. Soc. (3) 113 (2016) 627--648
\bibitem{mu}
M. Murugan. On the length of chains in a metric space. J. Func. Anal. 279 (2020) 108627
\end{thebibliography}
\end{document}